\documentclass{fundam}

\usepackage{url} % takes care of hyperlinks, preferred over hyperref
\usepackage[ruled,lined]{algorithm2e}% provides Algorithm environment
\usepackage{graphicx}% allows for inclusion of graphic files (figures)
\usepackage{tikz}
\usetikzlibrary{automata, positioning, arrows}
\usepackage{amssymb}
\usepackage{listings}

\newtheorem{Thm}{Theorem}
\newtheorem{Def}{Definition}
\newtheorem{Prob}{Problem}
\newtheorem{Cor}{Corollary}
\newtheorem{Rem}{Remark}

\newtheorem{Exm}{Example}

\begin{document}

%%%%%%%  parameters to be filled in by copy-editor  %%%%%%%%%%

\setcounter{page}{145}
\publyear{22}
\papernumber{2156}
\volume{189}
\issue{2}

\finalVersionForARXIV
%\finalVersionForIOS

%%%%%%%%%%%%%%%%%%%%%%%%%%%%%%%%%%%%%%

\title{On Taxicab Distance Mean Functions and their Geometric Applications: Methods, Implementations and Examples}

\address{Institute of Mathematics, University of Debrecen, P.O.Box 400, H-4002 Debrecen, Hungary}

\author{Csaba Vincze\thanks{Research supported in part by the E\"otv\"os Lor\'and Research Network (ELKH).}\\
Institute of Mathematics \\
University of Debrecen \\ P.O.Box 400, H-4002 Debrecen, Hungary\\
csvincze{@}science.unideb.hu
\and \'Abris Nagy \\
Institute of Mathematics \\
University of Debrecen \\
Debrecen, Hungary \\
abris.nagy{@}science.unideb.hu}

\maketitle

\runninghead{C. Vincze and \'A. Nagy}{On Taxicab Distance Mean Functions and their Geometric Applications}

\begin{abstract}
A distance mean function measures the average distance of points from the elements of a given set of points (focal set) in the space. The level sets of a distance mean function are called generalized conics. In case of infinite focal points the average distance is typically given by integration over the focal set. The paper contains a survey on the applications of taxicab distance mean functions and generalized conics' theory in geometric tomography: bisection of the focal set and reconstruction problems by coordinate X-rays. The theoretical results are illustrated by implementations in Maple, methods and examples as well. \footnote{The paper is based on the plenary lecture presented at Meeting on Tomography and Applications (Discrete Tomography, Neuroscience and Image Reconstruction) 16th Edition, IN MEMORIAM OF CARLA PERI, 2 - 4 May 2022, Mathematics Department, Politecnico di Milano, Milano, Italy.}
\end{abstract}


\section{Introduction}

A distance mean function measures the average distance of points from the elements of a given set of points (focal set) in the space. The level sets of a distance mean function are called generalized conics. The most important discrete examples are polyellipses (polyellipsoids) as level sets of a function measuring the arithmetic mean of distances from finitely many focal points (constant distance sum) and polynomial lemniscates as level sets of a function measuring the geometric mean of distances from finitely many focal points (constant distance pro\-duct). In case of infinite focal points the average distance is typically given by integration over the focal set. Using partitions and integral sums, the level sets (generalized conics) are Hausdorff  limits of polyellipsoids  Section 4 \cite{VinKovCsor}.

\subsection{Some general observations} The general form of the functions we are interested in is
\begin{equation}
\label{genformula}
x \mapsto f_D(x):=\frac{1}{\mu(D)}\int_{D} u\circ d(x,y)\, d_{\mu} y,
\end{equation}
where $d$ measures the distance between the points, $D\subset \mathbb{R}^n$ is a compact subset with a finite positive measure with respect to $\mu$, $u\colon \mathbb{R}\to \mathbb{R}$ is a strictly monotone increasing convex function satisfying the initial condition $u(0)=0$. In what follows we suppose that the distance function comes from a norm. The convexity of the integrand implies that the distance mean function is a convex and, consequently, a continuous function. Using the increasing slope property of convex functions, we have that
\begin{equation}
\label{ucoercive}
\liminf_{t\to \infty} \frac{u(t)}{t}>0
\end{equation}
and the distance mean function inherits a growth property of the form
\begin{equation}
\label{growth}
\liminf_{\|x\|\to \infty}\ \frac{f_D(x)}{\|x\|} > 0.
\end{equation}
The growth property (\ref{growth}) implies that the sublevel sets of the form $C_D:=\{ x \ | \ f_D(x)\leq  c\}\subset \mathbb{R}^n$ are bounded because the existence of a sequence $x_n\in C_D$ such that $\lim_{n\to \infty} \|x_n\|=\infty$ gives a contradiction:
$$\lim_{n\to \infty}\ \frac{f_D(x_n)}{\|x_n\|}\leq \lim_{n\to \infty}\ \frac{c}{\|x_n\|}=0.$$

\begin{Thm}{\cite{VinNagy12, VinNagy17}}
\label{comconv} The sublevel sets of a distance mean function are convex and compact.
\end{Thm}

Weierstrass theorem states that if all the level sets of a continuous function defined on a non-empty closed set in $\mathbb{R}^n$ are bounded, then it has a global minimizer.

\begin{Thm}{\cite{VinNagy12, VinNagy17}}
The distance mean function has a global minimizer.
\end{Thm}

\subsection{The problem of unicity} The general problem of unicity means to characterize the subsets in the space that are uniquely determined by the average distance measuring. The following theorem shows that the iteration of the averaging process determines the sublevel sets of a distance mean function under the choice $u(t)=t$ in formula (\ref{genformula}). The result is a generalization of \cite[Theorem 9]{VinNagy12}.

\begin{Thm} {\cite{VinNagy17}}
\label{unicity1} Let $u(t)=t$ and let $C_D$ be a sublevel set of $f_D$. If $f_{C_*}=f_{C_D}$, where $C_*$ is a compact set such that $\mu(C_D)=\mu(C_*)$ then $C_D$ is equal to $C_*$ except on a set of measure zero with respect to $\mu$.
\end{Thm}

We present some applications under a more special choice of the ingredients (\ref{genformula}). In all the following sections we choose $u(t)=t$ and $d=d_1$, the taxicab distance. In sections 2 and 3 we let $\mu$ be the Lebesgue measure, while $\mu$ is the counting measure in section 4.

\section{Bisection of bodies by coordinate hyperplanes} Suppose that distance measuring and integration are taken with respect to the taxicab distance
\begin{equation}
d_1(x,y)=\sum_{i=1}^n|x^i-y^i|
\end{equation}
and the Lebesgue measure $\mu_n$, respectively. Let $K$ be a compact subset of measure one\footnote{It is a technical condition to avoid the denominator $\mu_n(K)$.} in $\mathbb{R}^n$  and consider the taxicab distance mean function
\begin{equation}
\label{unweighted} f_{K}(x)=\int_K d_1(x,y)\, dy =\sum_{i=1}^n\int_{K} |x^i-y^i|\, d y.
\end{equation}
Since the derivative of the integrand at $x^i$ is $\pm 1$ depending on $y^i < x^i$ or $x^i< y^i$, we can conclude that the value $1$ occurs as many times as many points $y\in K$ is on the left hand side of $x$ with respect to the $i$-th coordinate:
$$K<_i x^i:=\{y\in K \ | \ y^i< x^i \}.$$
In a similar way, $-1$ occurs as many times as many points $y\in K$ is on the right hand side of $x$ with respect to the $i$-th coordinate:
$$x^i<_i K:=\{y\in K \ | \ x^i< y^i\}.$$
Since the set
$$K=_i x^i:=\{y\in K \ | \ y^i=x^i \} \ \ (i=1, \ldots, n)$$
is of measure zero we have that
\begin{equation}
\label{gradient}
D_i f_K(x)=\mu_{n}(K\leq_i x^i)-\mu_{n}(x^i\leq_i K)\ \ (i=1, \ldots, n).
\end{equation}

\begin{Thm}{\cite{VinNagy12}}
The point $x\in \mathbb{R}^n$ is a minimizer of $f_K$ if and only if each coordinate hyperplane at $x$ divides $K$ in two parts of equal measure.
\end{Thm}
\eject

How to bisect a set in two parts of equal measure? Formula (\ref{gradient}) shows that
$$|D_if_K(x)-D_if_K(y)|=2\mu_n \left(\min \{x^i, y^i\} <_i K <_i \max \{x^i, y^i\}\right)$$
and the compactness of $K$ implies that $f_K$ has a Lipschitzian gradient. Therefore the gradient descent method can be used to find the minimizer bisecting the measure of the integration domain $K$ in the sense that each coordinate hyperplane passing through the minimizer divides the set into two parts of equal measure. Let us present the gradient descent method in terms of a stochastic algorithm \cite{BNNV, VinNagy17}: let $P_k$ be a sequence of $K$-valued independent uniformly distributed random variables and consider the recursion
\begin{equation}
\label{random}
X_{k+1}=X_k-t_{k+1}Q_{k+1},
\end{equation}
where $X_0\in K$ is a (random) starting point,
\begin{equation}
\label{stepdirection}
Q_{k+1}:=\left(\textrm{sgn\ } (X_k^1-P_{k+1}^1), \ldots, \textrm{sgn\ } (X_k^n-P_{k+1}^n)\right)
\end{equation}
and the step size is a decreasing sequence of positive real numbers $t_k$ satisfying conditions
\begin{equation}
\label{stepsize}
\sum_{k=1}^{\infty}t_k=\infty\ \ \textrm{and}\ \ \sum_{k=1}^{\infty}t_k^2 < \infty.
\end{equation}
Assuming that $K$ is of measure one, we have the conditional probability
\begin{equation}
P(Q_{k+1}=(1, \ldots, 1)| X_k)=\mu_n \left((K < X_k^1) \cap \ldots \cap (K < X_k^n)\right)
\end{equation}
because  $Q_{k+1}=(1, \ldots, 1)$ means that $X_{k}$ is greater than $P_{k+1}$ with respect to the coordinatewise partial ordering $x\prec y \ \Leftrightarrow \ x^1 < y^1, \ldots, x^n < y^n$  and $P_{k+1}$ is a uniformly distributed $K$-valued random variable. In a similar way we have the conditional probability
\begin{equation}
P(Q_{k+1}=(1, -1, 1, \ldots, 1)| X_k)=
\end{equation}
$$\mu_n \left((K < X_k^1) \cap (X_k^2 < K)\cap (K < X_k^3)\cap \ldots \cap (K < X_k^n)\right), \ldots$$
and so on. A direct computation shows that $\displaystyle{\mathbb{E}(Q_{k+1}|X_k)=\textrm{grad\ } f_K(X_k)}.$

\begin{Rem}
To illustrate the process let us consider the case of dimension two. The lines parallel to the coordinate axis at $X_k$ divide the plane into four quadrants. Since we  have a sequence of independent, uniformly distributed random variables, the value of $P_{k+1}$ is most likely to fall in the quadrant containing the part of $K$ of the highest measure. Using formula (\ref{gradient}), it follows that the gradient of $f_K$ at $X_k$ is pointed in the same quadrant represented by the value of the stochastic vector $Q_{k+1}$. Therefore the step of the highest probability is taken into the opposite direction of the gradient in the sense that the corresponding quadrants are opposite to each other.
\end{Rem}

\begin{Def}
A nonempty compact set is called a body if it is the closure of its interior.
\end{Def}

\begin{Thm} {\cite{BNNV, VinNagy17}}
Let $K\subset \mathbb{R}^n$ be a connected compact body. The sequence of random variables $X_k$ converges almost surely to the unique global
minimizer $x^*$ of the function $f_K$.
\end{Thm}

\subsection{Implementation and examples}

We show an implementation of the above procedure in the Maple software for polygons in the plane. For this we first need to know how to choose a random point uniformly in a polygon. The standard way to do this is the following:
\begin{enumerate}
\itemsep=0.8pt
	\item Triangulate the polygon with the help of non-intersecting diagonals.\vspace*{-1mm}
	\item Compute the areas of the triangles and then the area of the whole polygon. \vspace*{-1mm}
	\item Choose a random number $x$ uniformly from the interval $[0,A]$, where $A$ denotes the area of the polygon. \vspace*{-1mm}
	\item If $A_i$ denotes the area of the triangle $T_i$ in the triangulation, then find the smallest positive integer $k$ that satisfies
		\[x\leq \sum_{i=1}^{k} A_i. \vspace*{-1mm}
	\]
	\item Then choose a random point uniformly in the triangle $T_k$. This can be done by choosing two independent random numbers $u$ and $v$ with uniform distribution in the $[0,1]$ interval. If the vertices of $T_k$ are denoted by $P_k$, $Q_k$, and $R_k$, and $u+v\leq 1$, then choose the point $P_k+u(Q_k-P_k)+v(R_k-P_k)$. If $u+v> 1$, then choose the point $P_k+(1-u)(Q_k-P_k)+(1-v)(R_k-P_k)$ of the triangle.
\end{enumerate}
Several algorithms exist for polygon triangulation. The most widely used algorithm is based on partitioning the polygon into monotone pieces first and then triangulating the monotone pieces \cite{ Garey, LeePrep}. Another favorable algorithm is the ear-clipping method based on the two ears theorem \cite{Meisters}. Both of them are built in the computational geometry package of Maple. Using this, the following procedure in Maple produces $n$ random points with uniform distribution in the polygon $P$ given by listing the vertices along the boundary.

\footnotesize
\begin{lstlisting}
> randompointsinpolygon:=proc(P,n)
local L,T,triangleareas,areasum,i,k,x,tf,u,v,Px,Py;
L:=[];
T:=ComputationalGeometry:-PolygonTriangulation(P);
triangleareas:=[];
for i from 1 to nops(T) do
  geometry:-point(A,P[T[i][1]][1],P[T[i][1]][2]);
  geometry:-point(B,P[T[i][2]][1],P[T[i][2]][2]);
  geometry:-point(C,P[T[i][3]][1],P[T[i][3]][2]);
  geometry:-triangle(t,[A,B,C]);
  triangleareas:=[op(triangleareas),geometry:-area(t)];
end do;
areasum:=add(evalf(triangleareas[i]),i=1..nops(triangleareas));
for j from 1 to n do
  x:=RandomTools:-Generate(float(range=0..areasum,method=uniform));
  k:=0;
  tf:=true;
  while tf do
    k:=k+1;
    if x<=add(evalf(triangleareas[i]),i=1..k) then
      tf:=false;
    end if;
  end do;
  u:=RandomTools:-Generate(float(range=0..1,method=uniform));
  v:=RandomTools:-Generate(float(range=0..1,method=uniform));
  if 1<u+v then
    u:=1-u;
    v:=1-v;
  end if;
  Px:=P[T[k][1]][1]+u*(P[T[k][2]][1]-P[T[k][1]][1])
	+v*(P[T[k][3]][1]-P[T[k][1]][1]);
  Py:=P[T[k][1]][2]+u*(P[T[k][2]][2]-P[T[k][1]][2])
	+v*(P[T[k][3]][2]-P[T[k][1]][2]);
  L:=[op(L),[Px,Py]];
end do;
return(L)
end proc;
\end{lstlisting}

\normalsize

\begin{figure}[!h]
\vspace*{-6mm}
\centering
\includegraphics[height=5.9cm]{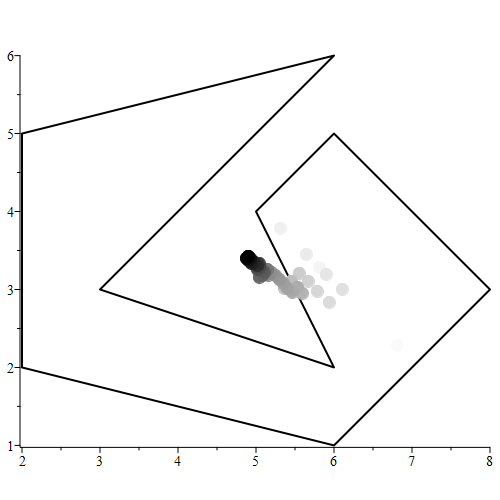}\quad\quad \includegraphics[height=6cm]{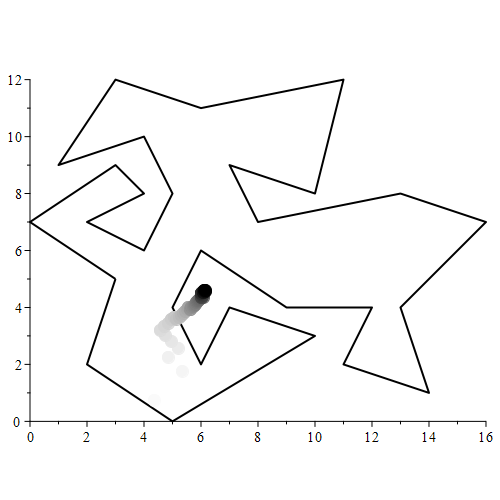}\vspace*{-4mm}
\caption{The sequence of points $X_k$ generated by the above procedure for $k=1,2,\ldots 50$. Darker points present elements $X_k$ with higher indices $k$. Notice how these sequences of points converge to the minimizer of the taxicab distance mean function \eqref{unweighted}.}
\end{figure}

Then the stochastic algorithm for finding the minimizer of the taxicab distance mean function \eqref{unweighted} of a polygon can be implemented in Maple as follows, see Figure 1. The step size sequence is given by $t_k=1/k$.

\footnotesize

\begin{lstlisting}
> minimizer:=proc(polygon,n)
local P,X,k,Q;
P:=randompointsinpolygon(polygon,n);
X:=P[1];
for k from 1 to nops(P)-1 do
  Q:=[signum(X[1]-P[k+1][1]),signum(X[2]-P[k+1][2])];
  X:=[X[1]-1/k*Q[1],X[2]-1/k*Q[2]];
end do;
return(X)
end proc;
\end{lstlisting}

\normalsize

\section{Applications in geometric tomography} The unweighted function \eqref{unweighted} is strongly related to the parallel X-rays as follows: by the Cavali\'{e}ri principle, formula
$$D_i f_K(x)=\mu_{n}(K\leq_i x^i)-\mu_{n}(x^i\leq_i K)\ \ (i=1, \ldots, n)$$
of the first partial derivatives implies that
\begin{equation}
\label{difformula}D_iD_if_K(x)=_{\textrm{a.e.}}2X_i K(x^i)\ \ (i=1, \ldots, n),
\end{equation}
where $X_iK(x^i):=\mu_{n-1}(x^i=_i K)$ is the $(n-1)$-dimensional Lebesgue measure of the set
\begin{equation}
\label{slice}
x^i=_i K:=\{y\in K\ | \ y^i=x^i\}.
\end{equation}
The functions
$$X_iK(t):=\mu_{n-1}(t=_i K) \ \ \ (t\in \mathbb{R} \ \textrm{and} \ i=1, \ldots, n)$$
are called the coordinate X-rays of $K$, see Figure 2. In terms of the coordinate X-rays
\begin{equation}
\label{intformula}
f_K(x)=\int_K d_1(x,y)\, dy =\sum_{i=1}^n\int_{K} |x^i-y^i|\, d y=\sum_{i=1}^n\ \intop\limits_{-\infty}^{\infty}|x^i-t|X_iK(t)\, dt.
\end{equation}
\begin{Thm} {\cite{VinNagy12}} $f_K=f_{L}$ iff the coordinate X-rays of $K$ and $L$ coincide almost everywhere.
\end{Thm}

\begin{figure}
\vspace*{2mm}
\centering
\includegraphics[height=5.3cm]{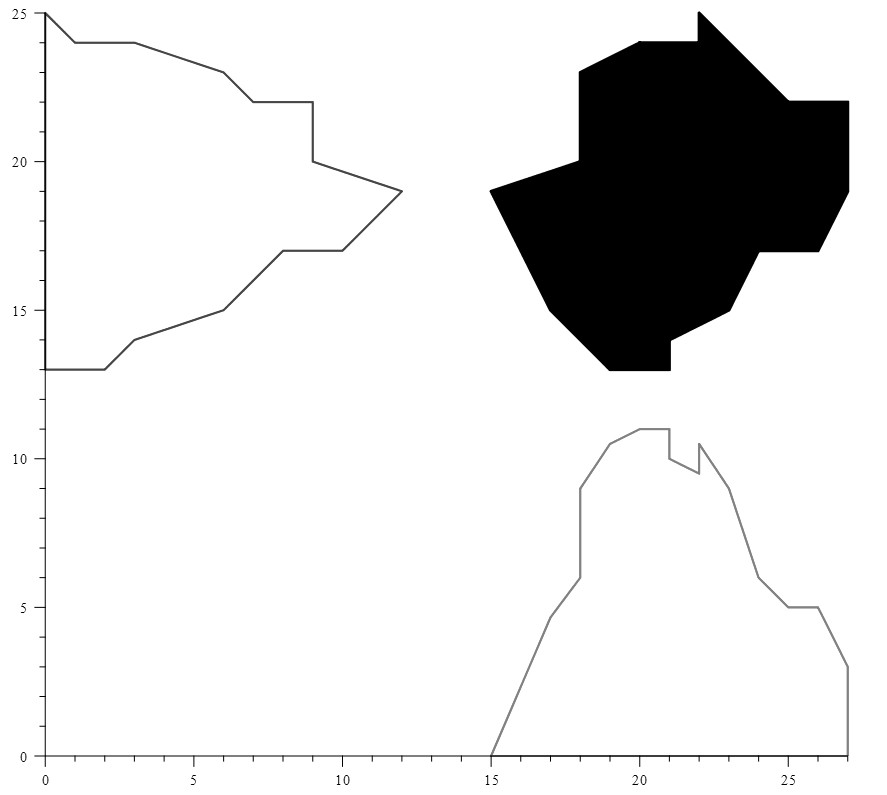}\quad\quad \includegraphics[height=5.3cm]{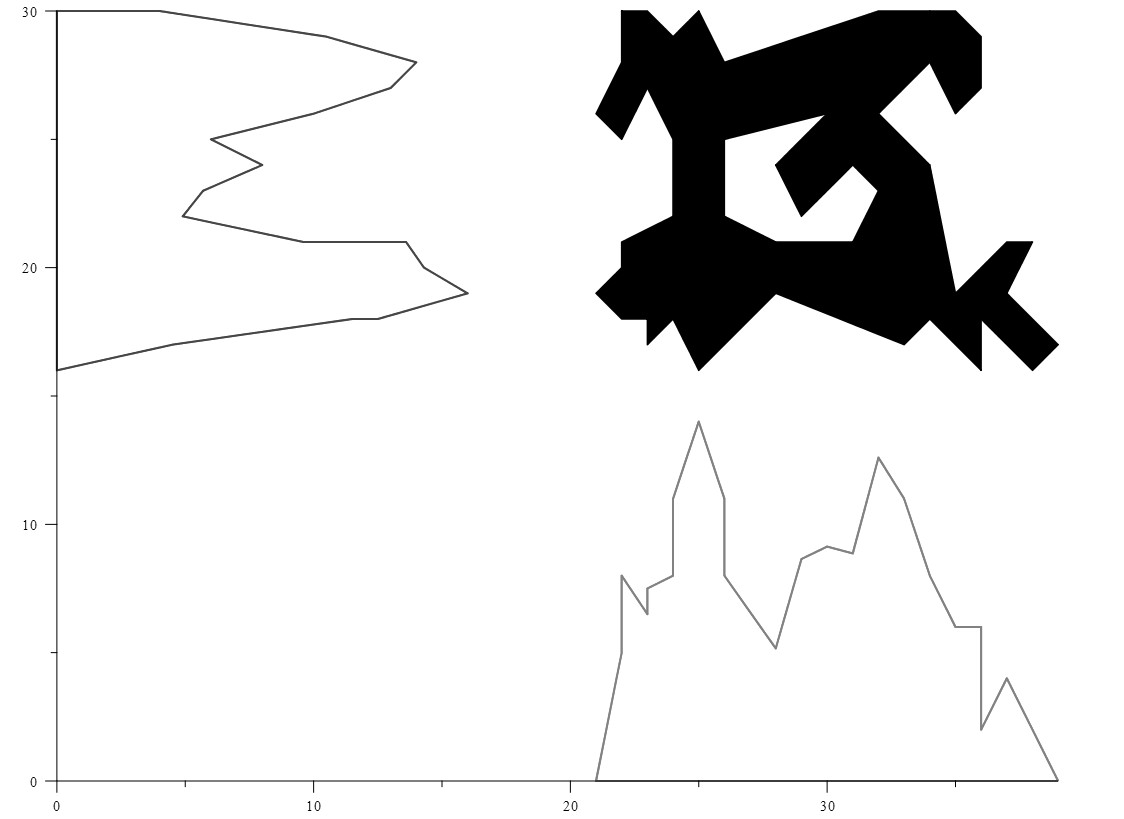}\vspace*{-1mm}
\caption{X-rays of compact planar bodies.}
\end{figure}

Since the coordinate X-rays determine both the measure and the taxicab distance mean function of the sets, we can formulate the following result as a consequence of Theorem \ref{unicity1}.

\begin{Cor} {\cite{VinNagy12, VinNagy17}}
The sublevel sets of a taxicab distance mean function are determined by their X-rays pa\-rallel to the coordinate hyperplanes among compact sets.
\end{Cor}

\begin{Exm} Circles are determined by their X-rays in the coordinate directions among compact sets in the plane. They are level sets of the taxicab distance mean function $f_B$ associated to the circumscribed square: if
$B:=\textrm{conv} \left \{(0,0), (1,0), (1,1), (0,1)\right\}$,
then we have that
$$f_B(x)=\big(x^1-(1/2)\big)^2+\big(x^2-(1/2)\big)^2+(1/2)$$
for any interior point $(x^1, x^2)\in B$. The general problem of unicity for convex bodies can be found in Gardner's basic monograph \cite{G06}: characterize those convex bodies that can be determined by two X-rays.
\end{Exm}

The taxicab distance mean function $f_K$ accumulates the coordinate X-ray information. Instead of the X-rays we can investigate a convex function independently of the convexity of the integration domain. Techniques and results based on $f_K$ are typically working in higher dimensional spaces as well.

\subsection{Reconstruction of planar sets by their coordinate X-rays} \cite{VinNagy12, VinNagy14a, VinNagy14} In what follows we restrict ourselves to the coordinate plane $\mathbb{R}^2$. Let $K\subset \mathbb{R}^2$ be a compact subset. The coordinate X-rays of $K$ enable us to construct an axis-parallel bounding box containing $K$. Since $f_K$ is also given by the coordinate X-rays, the reconstruction is based on the best approximation of $f_K$ by the distance mean functions of a special class of sets. They are constituted by unions of subrectangles of the bounding box under a given resolution: $f_{L_n}\to f_K.$ Taking $K^*$ as the limit set of a convergent subsequence in $L_n$, we have to provide the continuity of the mapping $L\mapsto f_L$ for a convergent reconstruction process. The continuity implies that the taxicab distance mean functions of $K^*$ and $K$ coincide. So do their coordinate X-rays (almost everywhere). In general X-rays can have deviant behavior under the Hausdorff convergence of the sets \cite{VinNagy16}. The taxicab distance mean functions are more regular objects in some sense. This makes them to be a natural starting point of the reconstruction.

\medskip
Let $K$ be a compact subset in the plane. The outer parallel body $K_{\varepsilon}$ is the union of closed Euclidean
disks centered at the points of $K$ with radius $\varepsilon > 0$. The Hausdorff distance between the compact subsets $K$ and $L$ is given by the formula
$$\delta (K,L) := \inf \{\varepsilon > 0 \ |\  K\subset L_{\varepsilon} \ \textrm{and}\ L\subset K_{\varepsilon}\}.$$

\begin{Def} \cite{VinNagy12} The Hausdorff convergence $L_n\to K$ is called regular if and only if
$$\lim_{n\to \infty}\mu_2(L_n)=\mu_2(K).$$
It is X-regular if and only if $\ \lim_{n\to \infty}\mu_2(I_n)=\mu_2(K)$, where $I_n:=\cap_{i=n}^{\infty} L_i$.
\end{Def}

It can be easily seen that under the hypothesis of the Hausdorff convergence, the regularity is equivalent to the convergence in the symmetric difference metric (or Lebesgue metric). In general they are not equivalent metrics as the following theorem shows.

\begin{Thm} {\cite{BBGV12}} The sequence $L_n$ converges in Hausdorff distance to $K$ if and only if
$$\lim_{n\to \infty} \mu_2 ((L_n)_{\varepsilon} \ \triangle \ K_{\varepsilon})=0\ \ \textrm{for each}\ \ \varepsilon>0,$$
where $K_{\varepsilon}$ is the parallel body of $K$ with radius $\varepsilon$.
\end{Thm}

\begin{Thm} {\cite{VinNagy12, VinNagy14}} If $L_n\to K$ with respect to the Hausdorff metric then
$$\limsup_{n\to \infty} f_{L_n}(x)\leq f_K(x).$$
If the Hausdorff convergence $L_n\to K$ is regular then $\lim_{n\to \infty} f_{L_n}( x)= f_K(x)$ and the convergence $f_{L_n}\to f_K$ is uniform over any compact subset in $\mathbb{R}^2$. If the Hausdorff convergence $L_n\to K$ is X-regular then it is regular and the coordinate X-rays converge to the coordinate X-rays of the limit set almost everywhere:
$$\lim_{n\to \infty} X_1L_n(t)=_{\textrm{a.e.}}X_1K(t), \ \lim_{n\to \infty} X_2L_n(t)=_{\textrm{a.e.}}X_2K(t).$$
\end{Thm}

We have the following examples:
\begin{itemize}
\item[(i)] If each $L_n$ is obtained from a compact set $L$ via finitely many Steiner symmetrizations and Euclidean isometries then the Hausdorff convergence $L_n\to K$ is regular \cite[Lemma 3.2]{BBGV12}.\item[(ii)] Any outer Hausdorff approximation $K\subset L_n \to K$ is X-regular \cite[Lemma 1, Remark 2]{VinNagy12}.
\item[(iii)]  Let $f\colon K\to D_1\subset \mathbb{C}$ be a homeomorphism, where $D_1$ denotes the unit disk of the plane centered at the origin. If the mapping $f$ is differentiable (in complex sense) at each inner point of $K$ then, by Mergelyan's theorem, $f$ can be approximated uniformly on $K$ by polynomials: $P_n\to f$. Therefore we have  an approximation of $K$ by polynomial lemniscate domains of the form $|P_n(z)|\leq c_n$ in the sense that the maximal connected components $L_n^*\subset K$ of the lemniscate domains tends to $K$ with respect to the Hausdorff metric. If $K$ has a boundary of measure zero then the Hausdorff convergece is X-regular \cite[Section 5]{VinNagy12}.
\item[(iv)] If $L_n$ is a sequence of compact connected hv-convex sets tending to the limit $K$ with respect to the Hausdorff metric, then the convergence is regular \cite[Section 3]{VinNagy14}.
\item[(v)] The Hausdorff convergence of compact convex subsets $L_n$ to $K$ with non-empty interior is X-regular \cite[Section 4.1]{VinNagy17}.
\end{itemize}

In the sense of the last example, the Hausdorff convergence in the class of compact convex sets (with non\-empty interior) implies the X-regularity and, by Theorem 8, the reconstruction can be based on direct comparisons of X-rays; see Gardner and Kiderlen \cite{GK07} (four directions, compact convex planar bodies). Indeed, if the sequence $L_n$ is constructed by the approximation of the X-rays of $K$, then the X-regularity implies that the X-rays of $L_n$ tend to the X-rays of the accumulation points which also equal to the X-rays of $K$ (almost everywhere). Example (iv) shows that the Hausdorff convergence in the class of compact connected hv-convex sets implies the regularity and the reconstruction can be based on direct comparisons of the taxicab distance mean functions.

\begin{Thm} \cite{VinNagy14} Let $\mathcal{M}_B^{hv}$ denote the set of non-empty compact connected hv-convex sets contained in the axis parallel bounding box $B\subset \mathbb{R}^2$ and let $K\in\mathcal{M}_B^{hv}$. For any $\varepsilon> 0$ there exists $\sigma>0$ such that whenever
$$\int_{B}|f_L(x)-f_K(x)|\, d{x}< \sigma$$
holds for $L\in\mathcal{M}_B^{hv}$, then there exists $K^*$, satisfying $\delta (L, K^*)< \varepsilon$ and $f_K=f_{K^*}$. Therefore $K$ and $K^*$ have the same coordinate X-rays almost everywhere.
\end{Thm}

\subsection{An algorithm for the reconstruction}

\cite{VinNagy14a}
Let $n\in\mathbb{N}$ be a natural number and suppose that the coordinate X-rays $X_1K$, $X_2K$ of a non-empty compact connected hv-convex planar body $K\subset\mathbb{R}^2$ are given. The Cartesian product of the supports of the coordinate X-rays gives a box
\begin{equation}
B=\textrm{supp}\ (X_1K)\times \textrm{supp}\ (X_2K)=[a,b]\times [c,d]
\end{equation}
containing $K$. The function $f_K$ associated to $K$ is defined by the formula
\begin{equation}
f_K(x)=\intop\limits_{-\infty}^{\infty}\left|x_1-t\right|X_1K(t)\,dt+\intop\limits_{-\infty}^{\infty}\left|x_2-t\right|X_2K(t)\,dt.
\end{equation}
Let
	\[t_i^1=a+i\frac{b-a}{n}\ \textrm{and}\ t_j^2=d-j\frac{d-c}{n}\quad (i, j=0,\ldots,n)
\]
be equally spaced points. The control grid $G^n_K:=\left\{y_{ij}\in B_K \ |\ i,j=1,\ldots,n\right\}$
consists of the centers of the subrectangles
\begin{equation}
\label{subrect}
B^n_{ij}=[t_{i-1}^1,t_i^1]\times [t_j^2,t_{j-1}^2],
\end{equation}
where $i,j=1,\ldots,n$. The feasible set $\mathcal{H}_n$ contains an element $L$ if and only if it is a compact connected hv-convex set which can be written as union of elements of the collection (\ref{subrect}) and
\begin{equation}
\label{eq:lin}
f_L(y_{ij})\geq f_K(y_{ij}) \ \ \textrm{for any} \ \ i,j=1,\ldots,n.
\end{equation}
For the output we choose $L_n\in \mathcal{H}_n$ that minimizes
\begin{equation}
\label{objective}
\sum_{i,j=1}^n\frac{f_{L_n}(y_{ij})-f_{K}(y_{ij})}{n^2},
\end{equation}
see Figures 3, 4 and 5. The procedure can be formulated in terms of a linear 0 - 1 programming because any element $L$ in the feasible set can be represented as a $0-1$ interval matrix by the variables $x_{kl}$ and $\overline{x}_{kl}=1-x_{kl}$, where $x_{kl}=1$ if
$B^n_{kl}\subset L$ and $x_{kl}=0$ otherwise ($k,l=1,\ldots,n$). The linearization of the constraints is based on \cite[chapters 11 and 12]{Lisun}. The applications of the  greedy or the antigreedy algorithmic paradigms are also possible \cite[sections 7 and 8]{VinNagy14a}. They are based on deleting the subrectangle which causes the extremal (the greatest or the least) average descent of $f_{L_n}$ at the control points. In general the antigreedy version increases the number of the possible outputs for making some voting processes more effective. The algorithm is adapted to finitely many and/or noisy measurements of the coordinate X-rays as well \cite{VinNagy16}.

\begin{figure}[!ht]
\begin{center}
		\includegraphics[scale=0.31]{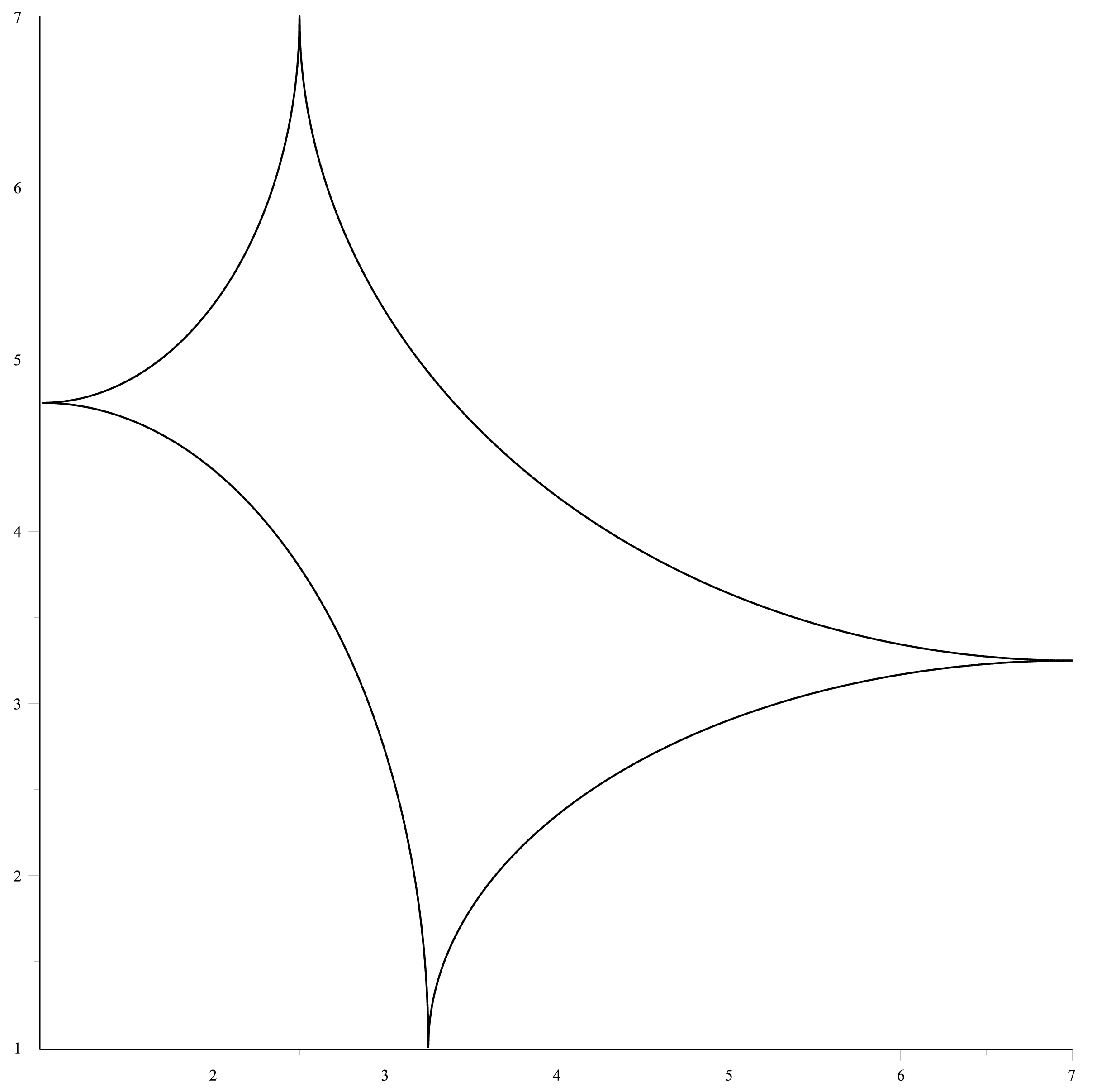}\vspace*{-1mm}
		\caption{The set we are looking for.}
\end{center}
\end{figure}

\begin{figure}[!h]
\vspace*{-3mm}
\begin{center}
\includegraphics[scale=0.32]{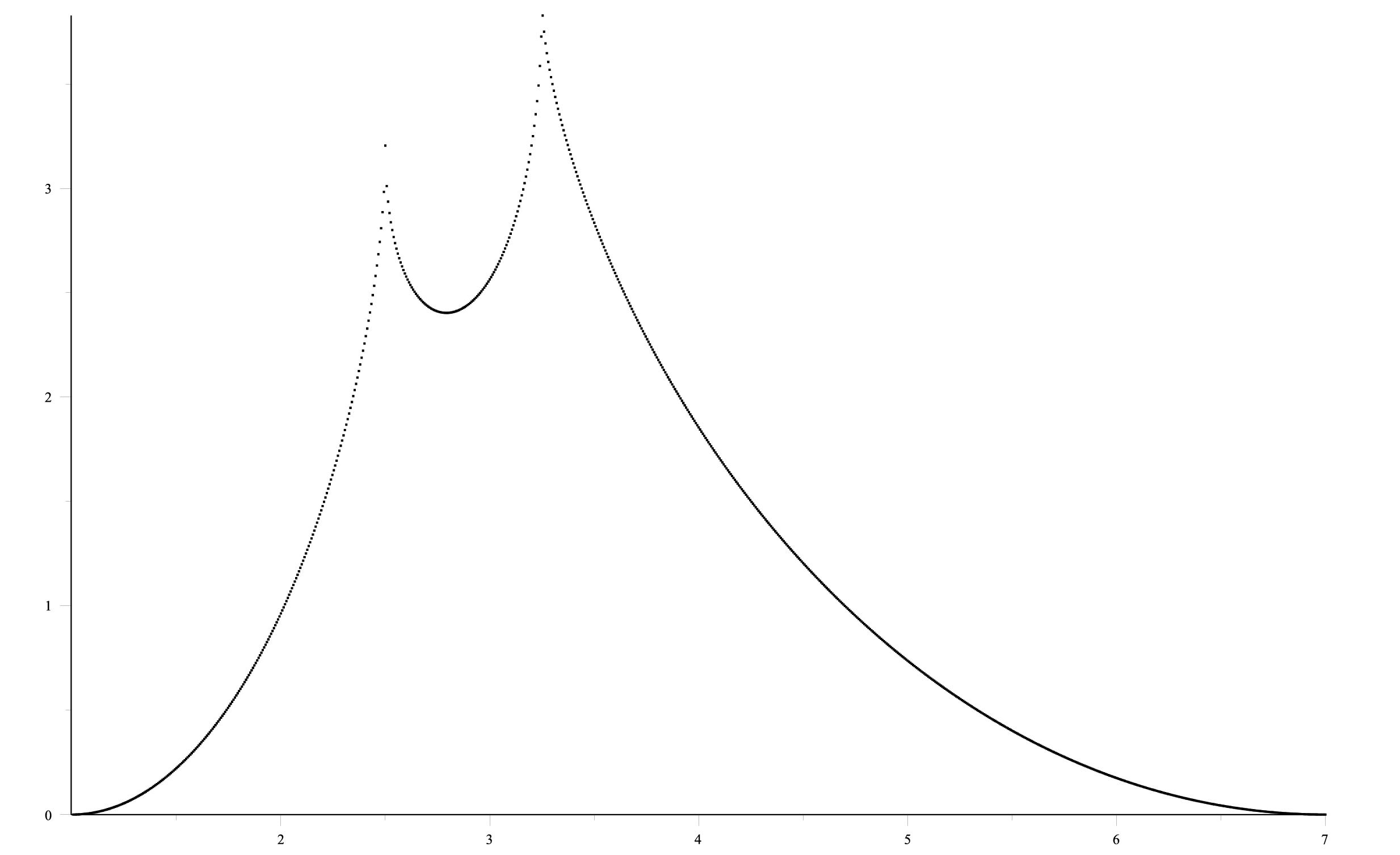} \quad \includegraphics[scale=0.32]{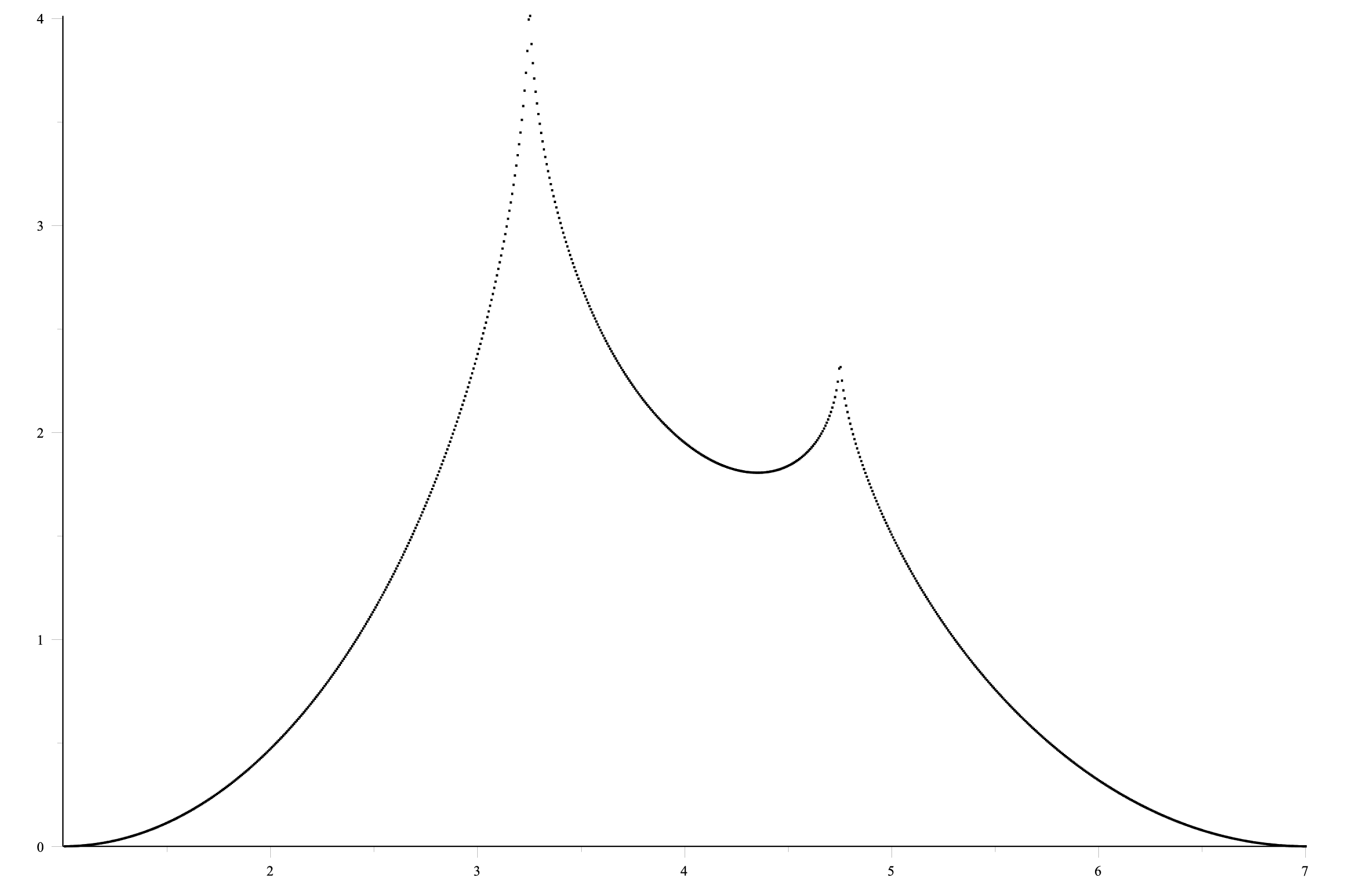}\vspace*{-1mm}
\caption{The coordinate X-rays.}
\end{center}
\end{figure}

\begin{figure}[!h]
\begin{center}
\includegraphics[scale=0.32]{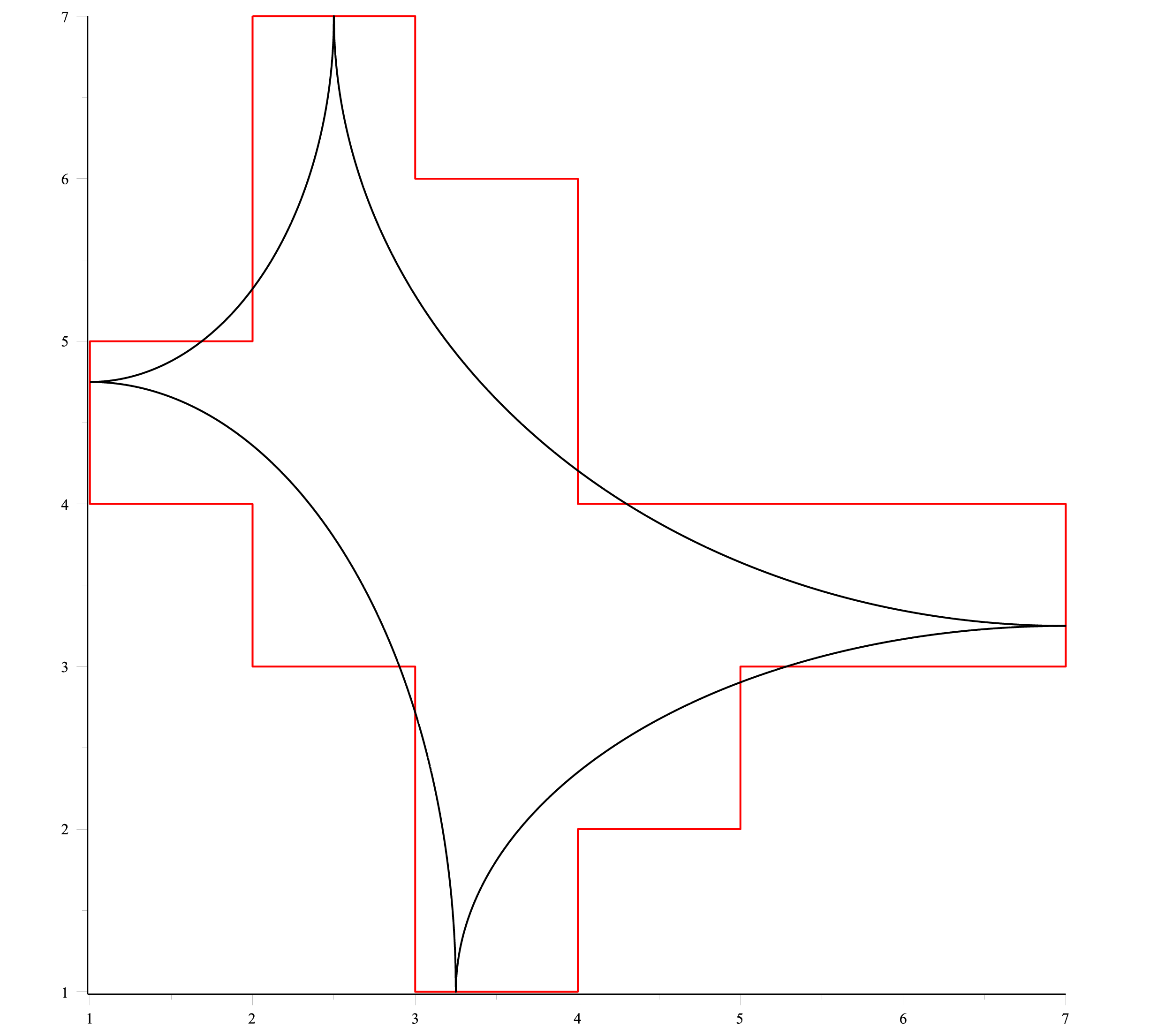} \quad \includegraphics[scale=0.32]{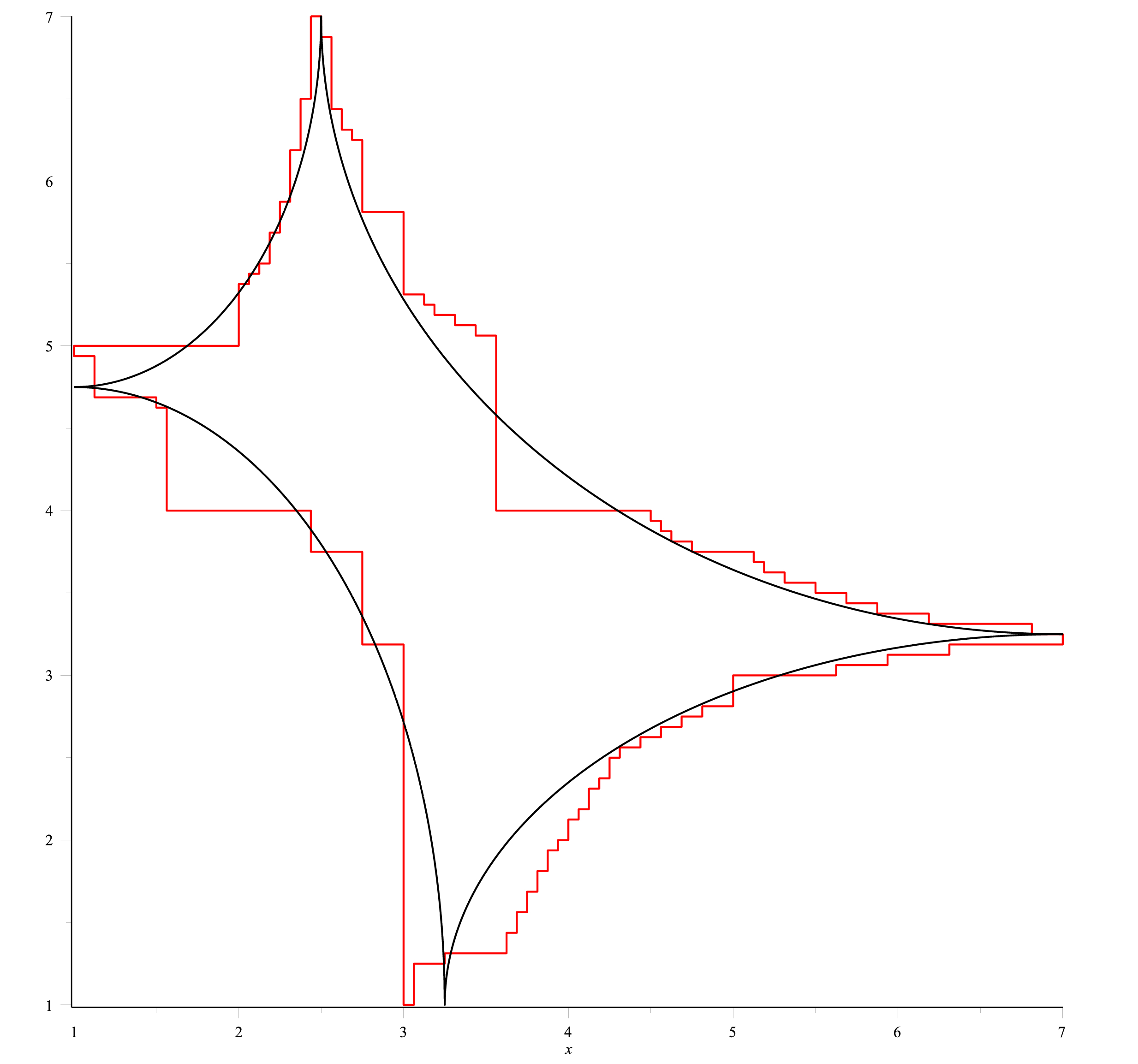}\vspace*{-1mm}
\caption{The optimal solution under low resolution (left) and a greedy version under high resolution (right).}
\end{center}
\end{figure}

\clearpage

\section{Reconstruction by the least average values} To illustrate the process let us consider the discrete version \cite{Vin} of the presented tomographic tools. It is a special case of the general theory with counting measure in the integral formulas. Let $F=\{x_i\in \mathbb{R}^n \ | \ i=1, \ldots, m\}$ be a finite set of different points in the coordinate space and consider the taxicab distance sum function
\begin{equation}
\label{sumfunction}
f(x):=\sum_{i=1}^m d_1(x, x_i)=\sum_{i=1}^m \sum_{j=1}^n|x^j-x_i^j|.
\end{equation}
Introducing the one-sided partial derivatives
$$D_j^+f(x):=\lim_{\varepsilon \to 0^+}\frac{f(x^1, \ldots, x^j+\varepsilon, \ldots, x^n)-f(x)}{\varepsilon},$$
$$D_j^-f(x):=\lim_{\varepsilon \to 0^-}\frac{f(x^1, \ldots, x^j+\varepsilon, \ldots, x^n)-f(x)}{\varepsilon}$$
we have the following collection of formulas:
$$ D_j^+f(x)=\left| F \leq _j x^j\right|-\left| F>_j x^j \right|,\ D_j^-f(x)=\left| F < _j x^j\right|-\left| F\geq _j x^j \right|,$$
where
$$F>_j t:=\{x_i\in F \ | \ x_i^j > t\}, \  F=_j t:=\{x_i\in F \ | \ x_i^j = t\}, \ F<_j t:=\{x_i\in F \ | \ x_i^j < t\},$$
$$F\geq _j t:=\{x_i\in F \ | \ x_i^j \geq t\}, \ F\leq _j t:=\{x_i\in F \ | \ x_i^j \leq t\},$$
$$\frac{D_j^+f(x)-D_j^-f(x)}{2}=\left | F=_j x^j \right | \quad (j=1, \ldots, n).$$
The cardinality $\left | F=_j x^j \right |$ is the number of the points in the intersection of $F$ with the hyperplane $x+H_j$, where $H_j:=\{x\in \mathbb{R}^n \ | \ x^j=0\}.$ The $(n-1)$-dimensional X-ray function parallel to the coordinate hyperplane $H_j$ is defined as
\begin{equation}
\label{xraydiscrete}
X_j\colon \mathbb{R}\to \mathbb{R}, \ X_j(t):=\left | F=_j t \right | \quad (j=1, \ldots, n).
\end{equation}
X-rays take the zero value except at finitely many $t\in \mathbb{R}$. In terms of X-rays
\begin{equation}\label{taxidistmean}
	f(x)=\sum_{i=1}^m d_1(x, x_i)=\sum_{i=1}^m \sum_{j=1}^n|x^j-x_i^j|=\sum_{j=1}^n \sum_{t\in \mathbb{R}} X_j(t)|x^j-t| \quad (x\in \mathbb{R}^n).
\end{equation}
Therefore the taxicab distance sum function accumulates the coordinate X-ray information.

\subsection{The least average value principle for the reconstruction of planar lattice sets by coordinate X-rays}

Let $G$ be the intersection of the integer lattice $\mathbb{Z}\times\mathbb{Z}$ and the rectangular picture region $[1,n]\times [1,m]$, where $m,n\in\mathbb{Z}$ are integers. Then we can write
	\[G=\{1, 2, \ldots, n\}\times \{1,2, \ldots, m\}\subset \mathbb{R}^2.
\]For each $j\in\left\{1,2,\ldots n\right\}$ let $l_{1,j}$ denote the vertical line intersecting the horizontal axis at the point $(j,0)$, and for each $i\in\left\{1,2,\ldots m\right\}$ let $l_{2,i}$ denote the horizontal line intersecting the vertical axis at the point $(0,m-i+1)$. The problem is to reconstruct the unknown lattice set $F\subset G$, if the numbers
	\[p_{1,j}=\left|l_{1,j}\cap F\right|,\quad j\in\left\{1,2,\ldots n\right\},
\]and
	\[p_{2,i}=\left|l_{2,i}\cap F\right|,\quad i\in\left\{1,2,\ldots m\right\}
\]are given. These are the numbers of elements of $F$ contained by each horizontal and vertical lattice line respectively, i.e. $p_{1,j}=X_1(j)$ and $p_{2,i}=X_2(m-i+1)$. The characteristic function of $F$ can be presented by the binary matrix $A=(a_{ij})$ of size $m\times n$, where
	\[a_{ij}=\begin{cases}
	1, &\textrm{if } (j,m-i+1)\in F,\\
	0, &\textrm{otherwise},
	\end{cases}
\]hence the above problem is equivalent to the following.

\begin{Prob}\label{tomography}
Given two integral vectors $R=(r_1,r_2,\ldots,r_m)$ and $S=(s_1,s_2,\ldots,s_n)$, find a binary matrix $A=(a_{ij})$ of size $m\times n$ such that\vspace*{-1mm}
	\[r_i=\sum_{j=1}^na_{ij},\quad i\in\left\{1,2,\ldots m\right\},
\]
and\vspace*{-1mm}
	\[s_j=\sum_{i=1}^ma_{ij},\quad j\in\left\{1,2,\ldots n\right\}.
\]
\end{Prob}
Certainly the above problem may have a solution only if $0\leq r_i\leq n$, $0\leq s_j\leq m$ for all $i\in\left\{1,2,\ldots,m\right\}$, $j\in\left\{1,2,\ldots,n\right\}$, and
	\[\sum_{i=1}^m r_i=\sum_{j=1}^n s_j.
\]If this property holds, then the integral vectors $R$ and $S$ are called compatible. In the rest of the paper we assume that the integral vectors $R$ and $S$ are compatible in Problem \ref{tomography}. This problem was first solved independently by Ryser \cite{Rys1, Rys2} and Gale \cite{Gale}. Ryser's approach is based on the construction of a maximal matrix\footnote{It is a matrix $(a_{ij})$ such that $a_{ij}=1$ whenever $j \leq r_i$, otherwise $a_{ij}=0$.} and shifting certain ones to the right within rows to attain the correct column sums. Gale's approach is based on network flows. This method was later improved by Anstee \cite{Anstee} and Batenburg \cite{Batenburg}. The advantage of the network flow method is that it can be easily applied for any pair of lattice directions, and even under the possible restriction, that we accept only those solutions, where a given set of entries are equal to zero and another given set of entries are equal to one. Before discussing the details of the network flow approach we introduce the basic concepts of flows in a network.

\medskip
Let $E$ and $V$ be two finite sets, and $\varphi\colon E\to V\times V$ a function. Then the triple $(E,V,\varphi)$ is called a directed graph, where $E$ is the set of edges, and $V$ is the set of vertices. If $\varphi(e)=(v_i,v_j)$ for some edge $e\in E$ and ordered pair of vetrices $(v_i,v_j)\in V\times V$, then we say $v_i$ is connected to $v_j$ by the edge $e$. The vertex $v_i$ is called the initial vertex, and $v_j$ is called the terminal vertex of the edge $e$. The directed path connecting the vertex $v_0$ to the vertex $v_k$ in a directed graph $(E,V,\varphi)$ is a sequence
	\[\left(v_0,e_1,v_1,e_2,v_2,e_3,v_3,\ldots v_{k-1},e_k,v_k\right),
\]where $v_0,v_1,\ldots v_k$ are pairwise different vertices, and $e_1,e_2,\ldots,e_k$ are pairwise different edges, such that $\varphi(e_i)=(v_{i-1},v_i)$ for all $i=1,2,\ldots k$. The undirected path connecting the vertex $v_0$ to the vertex $v_k$ in a directed graph $(E,V,\varphi)$ is similar to the directed path except, that now any of $\varphi(e_i)=(v_{i-1},v_i)$ or $\varphi(e_i)=(v_{i},v_{i-1})$ is possible for all $i=1,2,\ldots k$. An edge $e_i$ of an undirected path is called forward edge if $\varphi(e_i)=(v_{i-1},v_i)$, and it's called backward edge if $\varphi(e_i)=(v_{i},v_{i-1})$. We say that the length of a directed or undirected path is $k$ if it consists of $k$ edges. There are different efficient methods to find the shortest directed/undirected path connecting a vertex to another, for example with the help of breadth-first search.

\medskip
A network is a directed graph $(E,V,\varphi)$ together with a non-negative capacity function $U\colon E\to\mathbb{R}$ and two special vertices $s$ and $t$, such that there's no edge with terminal vertex $s$ and there's no edge with initial vertex $t$. Then the vertex $s$ is called source, while the vertex $t$ is called sink. The capacity of any edge $e\in E$ is denoted by $U(e)$. A flow on the network $(E,V,\varphi,U,s,t)$ is a function $Y\colon E\to \mathbb{R}$ which satisfies the following two conditions:
\begin{itemize}
	\item $0\leq Y(e)\leq U(e)$ for any edge $e$,
	\item for any vertex $v\in V$, except the source and the sink, it's true that the sum of the flow values on all the edges with terminal vertex $v$ is equal to the sum of the flow values on all the edges with initial vertex $v$.
\end{itemize}
The later is called the flow conservation property. The value of the flow on any edge $e\in E$ is denoted by $Y(e)$. We say that the edge $e$ is saturated if $Y(e)=U(e)$. The size of a flow $Y$ is the sum of the flow values on all the edges with the source $s$ being the initial vertex. The flow conservation property ensures that the size of the flow also equals to the sum of the flow values on all the edges with the sink $t$ being the terminal vertex.

\medskip
Now we discuss how a network is constructed for Problem \ref{tomography}, and how a flow of maximal size on the network helps to determine a solution of Problem \ref{tomography}. Let's define the network $(E,V,\varphi,U,s,t)$ in the following way:
\begin{itemize}
\itemsep=0.95pt
	\item a vertex $v_i$ is assigned to each horizontal line $l_{2,i}$,
	\item a vertex $w_j$ is assigned to each vertical line $l_{1,j}$,
	\item each vertex $v_i$ is connected to every vertex $w_j$,
	\item the source $s$ is connected to every vertex $v_i$,
	\item every vertex $w_j$ is connected to the sink $t$.
\end{itemize}
Hence the vertex set is
	\[V=(s,t,v_1,v_2,\ldots,v_m,w_1,w_2,\ldots w_n),
\]the edge set is
	\[E=\left\{e_{ij}\,|\,i\in\left\{1,2,\ldots m\right\},\,j\in\left\{1,2,\ldots n\right\}\right\}\cup\left\{e_{si}\,|\,i\in\left\{1,2,\ldots m\right\}\right\}\cup\left\{e_{jt}\,|\,j\in\left\{1,2,\ldots n\right\}\right\},
\]where
	\[\varphi(e_{ij})=(v_i,w_j),\quad \varphi(e_{si})=(s,v_i),\quad \varphi(e_{jt})=(w_j,t).
\]Then the capacity function is defined as
\begin{itemize}
\itemsep=0.95pt
	\item $U(e_{ij})=1$, for all $i\in\left\{1,2,\ldots m\right\}$ and $j\in\left\{1,2,\ldots n\right\}$,
	\item $U(e_{si})=r_i$, for all $i\in\left\{1,2,\ldots m\right\}$,
	\item $U(e_{jt})=s_j$, for all $j\in\left\{1,2,\ldots n\right\}$.
\end{itemize}
If $Y$ is an integer flow (i.e. flow with integer values), then the flow values $Y(e_{ij})$ are all equal to 0 or 1, since the capacities of the edges $e_{ij}$ are all equal to 1. Thus, having an integer flow $Y$, we can construct the binary matrix $A=(a_{ij})$ of size $m\times n$ as
	\[a_{ij}=Y(e_{ij}),\quad \textrm{for all } i\in\left\{1,2,\ldots m\right\},\ j\in\left\{1,2,\ldots n\right\}.
\]Actually, there's a one-to-one correspondence between integer flows on the network and binary matrices of size $m\times n$. The flow conservation property shows that the binary matrix $A$ corresponding to the integer flow $Y$ has row sums equal to $Y(e_{si})$ and column sums equal to $Y(e_{jt})$, where $Y(e_{si})\leq r_i$ and $Y(e_{jt})\leq s_j$. Furthermore it's easy to see that no flow can have size larger than
	\[\sum_{i=1}^m U(e_{si})=\sum_{j=1}^n U(e_{jt}),
\]and whenever the size equals to the above number, then all the edges $e_{si}$ and $e_{jt}$ are saturated. Now our task is to find a flow $Y^*$ of maximal size on the network. All the capacities of the network are integer numbers, and it can be proved that the maximal flow on such network also has integer values. Hence we have the following theorem.

\begin{Thm}
Problem \ref{tomography} has a solution if and only if the size of the maximal flow $Y^*$ on the network $(E,V,\varphi,U,s,t)$ is equal to
	\[\sum_{i=1}^mr_i=\sum_{j=1}^ns_j.
\]Then the binary matrix $A=(a_{ij})$ of size $m\times n$ with
	\[a_{ij}=Y(e_{ij}),\quad \textrm{for all } i\in\left\{1,2,\ldots m\right\},\ j\in\left\{1,2,\ldots n\right\}
\]is a solution of Problem \ref{tomography}.
\end{Thm}

The construction of the maximal flow has two stages. In the first stage we construct an initial flow, and in the second stage we increase the size of the flow by changing the flow values along so called flow-augmenting paths. A flow-augmenting path for the flow $Y$ on the network $(E,V,\varphi,U,s,t)$ is an undirected path connecting the source $s$ to the sink $t$, such that the forward edges are all unsaturated and the values of the flow $Y$ on the backward edges are strictly positive. If all the capacities are integer numbers and $Y$ is an integer flow, then it's easy to see that increasing the flow values by 1 along the forward edges and decreasing the flow values by 1 along the backward edges results a new flow $\widehat{Y}$ with larger size. It's possible to prove that a flow on a network is maximal if and only if there exists no flow-augmenting path. Hence starting with an initial flow, such as the zero flow, we can find a maximal flow by searching for flow-augmenting paths and increasing the size of the flow as long as such path exists. The process is accelerated much, if we always choose a shortest flow-augmenting path, which can be efficiently found with help of breadth-first search.

The construction of the network associated to Problem \ref{tomography} implies that any flow-augmenting path consists of at least 3 edges. This shows that if we choose the zero flow as initial flow at the first stage, then we first try to increase the flow values along flow-augmenting paths of 3 edges by 1 as long as possible without exceeding the capacities. This means for the matrix $A$ that we first let $A$ be the zero matrix, and try to switch entries of $A$ to 1 as long as possible without exceeding the corresponding row and column sums. How further (i.e. longer) flow-augmenting paths can be found and what they imply on the matrix will be discussed in section 4.3. Now we mention that it's also an important question which ordering of the entries of $A$ is considered when we try to switch them to $1$. It turns out that putting different preferences on the entries has a large effect on the number of further flow-augmenting paths required later to attain the maximal flow.

\medskip
The preference on the entries of $A$ can be based directly on the row sums or column sums, but since the distance mean function has a nice connection to the coordinate X-rays, we can choose preferences upon the values of discrete version of the distance mean function, which is called taxicab distance sum function defined by formula (\ref{taxidistmean}). The value of the taxicab distance sum function corresponding to the unknown set $F$ at any point $x=(x^1,x^2)\in\mathbb{R}^2$ can be computed as
\begin{equation}
\begin{split}
f(x) & = \sum_{j=1}^{n} X_1(j)|x^1-j|+\sum_{i=1}^{m} X_2(i)|x^2-i| \\
 & =\sum_{j=1}^{n} s_j\cdot |x^1-j|+\sum_{i=1}^{m} r_{m-i+1}\cdot|x^2-i| \\
 & =\sum_{j=1}^{n} s_j\cdot|x^1-j|+\sum_{i=1}^{m} r_{i}\cdot|x^2-(m-i+1)|.
\end{split}
\end{equation}
The least average value principle means that points of the set $G$, or equivalently entries of the matrix $A$, with lowest taxicab distance sum values are preferred first. This implies the following steps:
\begin{enumerate}
\itemsep=0.9pt
	\item Create the preference list $L$ of entries of $A$ by sorting the entries into increasing order with respect to taxicab distance sum values at the corresponding points of $G$.
	\item Try to switch the entry of $A$ that corresponds to the first element of the preference list. This is possible unless the prescribed row sum or column sum corresponding to that element is zero. Then we say that the first element of the preference list is visited.
	\item After having a preference list with some visited elements, look for the first unvisited element of the list and try to switch the corresponding element of $A$. This is possible if none of the prescribed row sums and column sums are exceeded by the switch. We repeat this as long as the preference list has unvisited elements.
\end{enumerate}
Two accelerating methods can be applied here. The first one is that once we attain a row where the row sum of $A$ equals to the prescribed row sum, then we set all the unvisited entries of $A$ in that row to be visited without switching them to 1. These entries of $A$ remain equal to zero until all elements of the preference list are visited. A similar step can be applied for columns as well. The second accelerating technique is that once we attain a row where the difference of the prescribed row sum and the actual row sum of $A$ equals to the number of unvisited entries in that row, then we switch these unvisited entries of $A$ to 1 in those columns, where the prescribed column sum is not exceeded by the switch. Then we also set all the unvisited entries of $A$ in that row to be visited. A similar step can be applied for columns as well. Certainly applying these accelerating methods has a computational cost, but in change it reduces the number of further necessary switches along flow-augmenting paths.

\subsection{Example.} Now we present an example of Problem 1, and starting with the zero matrix $A$, we show how entries of $A$ are replaced by ones based on the preference list determined by the taxicab distance sum values. We also use the two accelerating method described above. Let $m=n=5$ and hence $G=\{1,2,3,4,5\}\times \{1,2,3,4,5\}.$ Furthermore, let the row sums and column sums be
	\[(r_1,r_2,r_3,r_4,r_5)=(3,1,4,4,2)\quad\textrm{and}\quad (s_1,s_2,s_3,s_4,s_5)=(4,3,1,4,2)
\]respectively. The matrix $A$ initially equals to the zero matrix. The matrix
	\[M:=\left[\begin{array}{ccccc}
54 & 48 & 48 & 50 & 60 \\[0.2cm]
46 & 40 & 40 & 42 & 52 \\[0.2cm]
40 & 34 & 34 & 36 & 46 \\[0.2cm]
42 & 36 & 36 & 38 & 48 \\[0.2cm]
52 & 46 & 46 & 48 & 58
\end{array}\right]
\]shows the values of the taxicab distance sum function $f$ at the points of the set $G$, where
	\[m_{kl}=f(l,m-k+1)=\sum_{j=1}^{n} s_j\cdot|l-j|+\sum_{i=1}^{m} r_i\cdot|k-i|
\]for all $k,l\in \{1,2,3,4,5\}$. The lowest taxicab distance sum function value is $34$. The first element of the preference list $L$ is the entry with the lowest taxicab distance sum function value. There are two such entries: $(3,2)$ and $(3,3)$. We choose, for example, the former one and switch $a_{32}$ to $1$ in the matrix $A$. Then we say that the entry $(3,2)$ is visited, and the next unvisited element of the preference list is $(3,3)$. We switch $a_{33}$ to $1$.
	\[\left[\begin{array}{ccccc}
0 & 0 & 0 & 0 & 0 \\[0.2cm]
0 & 0 & 0 & 0 & 0 \\[0.2cm]
0 & 0 & 0 & 0 & 0 \\[0.2cm]
0 & 0 & 0 & 0 & 0 \\[0.2cm]
0 & 0 & 0 & 0 & 0
\end{array}\right] \longrightarrow
\left[\begin{array}{ccccc}
0 & 0 & 0 & 0 & 0 \\[0.2cm]
0 & 0 & 0 & 0 & 0 \\[0.2cm]
0 & {\bf 1} & 0 & 0 & 0 \\[0.2cm]
0 & 0 & 0 & 0 & 0 \\[0.2cm]
0 & 0 & 0 & 0 & 0
\end{array}\right] \longrightarrow
\left[\begin{array}{ccccc}
0 & 0 & 0 & 0 & 0 \\[0.2cm]
0 & 0 & 0 & 0 & 0 \\[0.2cm]
0 & {\bf 1} & {\bf 1} & 0 & 0 \\[0.2cm]
0 & 0 & 0 & 0 & 0 \\[0.2cm]
0 & 0 & 0 & 0 & 0
\end{array}\right]
\]We see that the 3rd column sum of $A$ equals to $s_3=1$, thus we say that all entries in the 3rd column are visited, but we don't change them. Then we see that the difference of $r_4$ and the 4th row sum of $A$ equals to the number of unvisited entries in the 4th row, thus we switch all the unvisited entries to 1 in the 4th row and say that they are visited.
	\[\left[\begin{array}{ccccc}
0 & 0 & 0 & 0 & 0 \\[0.2cm]
0 & 0 & 0 & 0 & 0 \\[0.2cm]
0 & {\bf 1} & {\bf 1} & 0 & 0 \\[0.2cm]
0 & 0 & 0 & 0 & 0 \\[0.2cm]
0 & 0 & 0 & 0 & 0
\end{array}\right]\longrightarrow
\left[\begin{array}{ccccc}
0 & 0 & {\bf 0} & 0 & 0 \\[0.2cm]
0 & 0 & {\bf 0} & 0 & 0 \\[0.2cm]
0 & {\bf 1} & {\bf 1} & 0 & 0 \\[0.2cm]
0 & 0 & {\bf 0} & 0 & 0 \\[0.2cm]
0 & 0 & {\bf 0} & 0 & 0
\end{array}\right]
\longrightarrow
\left[\begin{array}{ccccc}
0 & 0 & {\bf 0} & 0 & 0 \\[0.2cm]
0 & 0 & {\bf 0} & 0 & 0 \\[0.2cm]
0 & {\bf 1} & {\bf 1} & 0 & 0 \\[0.2cm]
{\bf 1} & {\bf 1} & {\bf 0} & {\bf 1} & {\bf 1} \\[0.2cm]
0 & 0 & {\bf 0} & 0 & 0
\end{array}\right]
\]The next two unvisited elements of the preference list (i.e. unvisited elements with lowest taxicab distance sum function value) are $(3,4)$ and then $(2,2)$ where we switch the matrix to 1, and we call these entries visited.
	\[\left[\begin{array}{ccccc}
0 & 0 & {\bf 0} & 0 & 0 \\[0.2cm]
0 & 0 & {\bf 0} & 0 & 0 \\[0.2cm]
0 & {\bf 1} & {\bf 1} & 0 & 0 \\[0.2cm]
{\bf 1} & {\bf 1} & {\bf 0} & {\bf 1} & {\bf 1} \\[0.2cm]
0 & 0 & {\bf 0} & 0 & 0
\end{array}\right] \longrightarrow
\left[\begin{array}{ccccc}
0 & 0 & {\bf 0} & 0 & 0 \\[0.2cm]
0 & 0 & {\bf 0} & 0 & 0 \\[0.2cm]
0 & {\bf 1} & {\bf 1} & {\bf 1} & 0 \\[0.2cm]
{\bf 1} & {\bf 1} & {\bf 0} & {\bf 1} & {\bf 1} \\[0.2cm]
0 & 0 & {\bf 0} & 0 & 0
\end{array}\right] \longrightarrow
\left[\begin{array}{ccccc}
0 & 0 & {\bf 0} & 0 & 0 \\[0.2cm]
0 & {\bf 1} & {\bf 0} & 0 & 0 \\[0.2cm]
0 & {\bf 1} & {\bf 1} & {\bf 1} & 0 \\[0.2cm]
{\bf 1} & {\bf 1} & {\bf 0} & {\bf 1} & {\bf 1} \\[0.2cm]
0 & 0 & {\bf 0} & 0 & 0
\end{array}\right]
\]
We see that the 2nd row sum of $A$ equals to $r_2=1$ and the 2nd column sum of $A$ equals to $s_2=3$, thus we say that all entries in the 2nd row and 2nd column are visited, but we don't change them.
\[\left[\begin{array}{ccccc}
0 & 0 & {\bf 0} & 0 & 0 \\[0.2cm]
0 & {\bf 1} & {\bf 0} & 0 & 0 \\[0.2cm]
0 & {\bf 1} & {\bf 1} & {\bf 1} & 0 \\[0.2cm]
{\bf 1} & {\bf 1} & {\bf 0} & {\bf 1} & {\bf 1} \\[0.2cm]
0 & 0 & {\bf 0} & 0 & 0
\end{array}\right] \longrightarrow
\left[\begin{array}{ccccc}
0 & 0 & {\bf 0} & 0 & 0 \\[0.2cm]
{\bf 0} & {\bf 1} & {\bf 0} & {\bf 0} & {\bf 0} \\[0.2cm]
0 & {\bf 1} & {\bf 1} & {\bf 1} & 0 \\[0.2cm]
{\bf 1} & {\bf 1} & {\bf 0} & {\bf 1} & {\bf 1} \\[0.2cm]
0 & 0 & {\bf 0} & 0 & 0
\end{array}\right] \longrightarrow
\left[\begin{array}{ccccc}
0 & {\bf 0} & {\bf 0} & 0 & 0 \\[0.2cm]
{\bf 0} & {\bf 1} & {\bf 0} & {\bf 0} & {\bf 0} \\[0.2cm]
0 & {\bf 1} & {\bf 1} & {\bf 1} & 0 \\[0.2cm]
{\bf 1} & {\bf 1} & {\bf 0} & {\bf 1} & {\bf 1} \\[0.2cm]
0 & {\bf 0} & {\bf 0} & 0 & 0
\end{array}\right]
\] Now the difference of $r_1$ and the 1st row sum of $A$ equals to the number of unvisited entries in the 1st row, thus we switch all the unvisited entries to 1 in the 1st row and say that they are visited.
	\[ \left[\begin{array}{ccccc}
0 & {\bf 0} & {\bf 0} & 0 & 0 \\[0.2cm]
{\bf 0} & {\bf 1} & {\bf 0} & {\bf 0} & {\bf 0} \\[0.2cm]
0 & {\bf 1} & {\bf 1} & {\bf 1} & 0 \\[0.2cm]
{\bf 1} & {\bf 1} & {\bf 0} & {\bf 1} & {\bf 1} \\[0.2cm]
0 & {\bf 0} & {\bf 0} & 0 & 0
\end{array}\right] \longrightarrow
\left[\begin{array}{ccccc}
{\bf 1} & {\bf 0} & {\bf 0} & {\bf 1} & {\bf 1} \\[0.2cm]
{\bf 0} & {\bf 1} & {\bf 0} & {\bf 0} & {\bf 0} \\[0.2cm]
0 & {\bf 1} & {\bf 1} & {\bf 1} & 0 \\[0.2cm]
{\bf 1} & {\bf 1} & {\bf 0} & {\bf 1} & {\bf 1} \\[0.2cm]
0 & {\bf 0} & {\bf 0} & 0 & 0
\end{array}\right]
\]At this point we see that the 5th column sum of $A$ equals to the column sum $s_5$, while in the 1st, and 4th columns the differences of the prescribed column sums $s_1$, $s_4$ and the actual column sums of $A$ equal to the number of unvisited entries in the corresponding columns. Thus, we first say that all entries in the 5th column are visited, but we don't change them. Then we switch the unvisited entries $a_{31}$ and $a_{51}$ of the 1st column to 1 and the unvisited entry $a_{54}$ of the 4th column to $1$.
	\[\left[\begin{array}{ccccc}
{\bf 1} & {\bf 0} & {\bf 0} & {\bf 1} & {\bf 1} \\[0.2cm]
{\bf 0} & {\bf 1} & {\bf 0} & {\bf 0} & {\bf 0} \\[0.2cm]
0 & {\bf 1} & {\bf 1} & {\bf 1} & 0 \\[0.2cm]
{\bf 1} & {\bf 1} & {\bf 0} & {\bf 1} & {\bf 1} \\[0.2cm]
0 & {\bf 0} & {\bf 0} & 0 & 0
\end{array}\right] \longrightarrow
\left[\begin{array}{ccccc}
{\bf 1} & {\bf 0} & {\bf 0} & {\bf 1} & {\bf 1} \\[0.2cm]
{\bf 0} & {\bf 1} & {\bf 0} & {\bf 0} & {\bf 0} \\[0.2cm]
0 & {\bf 1} & {\bf 1} & {\bf 1} & {\bf 0} \\[0.2cm]
{\bf 1} & {\bf 1} & {\bf 0} & {\bf 1} & {\bf 1} \\[0.2cm]
0 & {\bf 0} & {\bf 0} & 0 & {\bf 0}
\end{array}\right] \longrightarrow
 \left[\begin{array}{ccccc}
{\bf 1} & {\bf 0} & {\bf 0} & {\bf 1} & {\bf 1} \\[0.2cm]
{\bf 0} & {\bf 1} & {\bf 0} & {\bf 0} & {\bf 0} \\[0.2cm]
{\bf 1} & {\bf 1} & {\bf 1} & {\bf 1} & {\bf 0} \\[0.2cm]
{\bf 1} & {\bf 1} & {\bf 0} & {\bf 1} & {\bf 1} \\[0.2cm]
{\bf 1} & {\bf 0} & {\bf 0} & {\bf 1} & {\bf 0}
\end{array}\right]
\]Finally all elements of the preference list are visited. The final matrix $A$ has row sums $(3,1,4,4,2)=(r_1,r_2,r_3,r_4,r_5)$ and has column sums $(4,3,1,4,2)=(s_1,s_2,s_3,s_4,s_5)$, hence no further augmentation is required. In other examples, especially for larger matrices, it's possible that further augmentation is required as some of prescribed row sums and column sums are not attained yet. Anyway, the least average value principle transforms the coordinate X-rays into some geometric information by the values of the taxicab distance sum function. They are working as probability-like quantities whenever the subsequent step of the algorithm is not determined by the X-rays.

\subsection{Switching chains}
Here we discuss how to find a flow-augmenting path for a flow $Y$ in the network $(E,V,\varphi,U,s,t)$ constructed for Problem 1, if $Y$ is not maximal, but there exists no flow augmenting path of length 3. Then there exists at least one unsaturated edge $e_{si}$ (connecting the source $s$ to the vertex $v_i$) and at least one unsaturated edge $e_{jt}$ (connecting the vertex $w_j$ to the sink $t$). However the edge $e_{ij}$ must be saturated for any such pair of unsaturated edges $e_{si}$ and $e_{jt}$, because there's no flow augmenting path of length 3. This means that the $i$-th row sum is less than $r_i$ and the $j$-th column sum is less than $s_j$, while $a_{ij}=1$ in the binary matrix $A$ corresponding to $Y$. Finding a (shortest) flow-augmenting path containing the edges $e_{si}$ and $e_{jt}$ is equivalent to finding a (shortest) switching chain, i.e. sequence of pair of indexes $(i_0,j_0),(i_1,j_1),\ldots (i_l,j_l)$ that satisfies the following conditions:
\begin{itemize}
\itemsep=0.9pt
  \item $l$ is an odd number,
	\item $i_0=i$ and $j_0=j$,
	\item $1\leq i_k\leq m$ and $1\leq j_k\leq n$ for all $k\in \left\{0,1,\ldots,l\right\}$,
	\item $i_k=i_{k+1}$ and $j_k\neq j_{k+1}$ for all even numbers $k\in \left\{0,1,\ldots,l-1\right\}$,
	\item $j_k=j_{k+1}$ and $i_k\neq i_{k+1}$ for all odd numbers $k\in \left\{1,2,\ldots,l-1\right\}$,
	\item $j_l=j_0$ and $i_l\neq i_0$,
	\item $a_{i_kj_k}=1$ for all even numbers $k\in \left\{1,2,\ldots,l\right\}$,
	\item $a_{i_kj_k}=0$ for all odd numbers $k\in \left\{1,2,\ldots,l\right\}$,
\end{itemize}

\begin{figure}[!b]
\vspace*{-3mm}
\centering
\includegraphics[width=8cm]{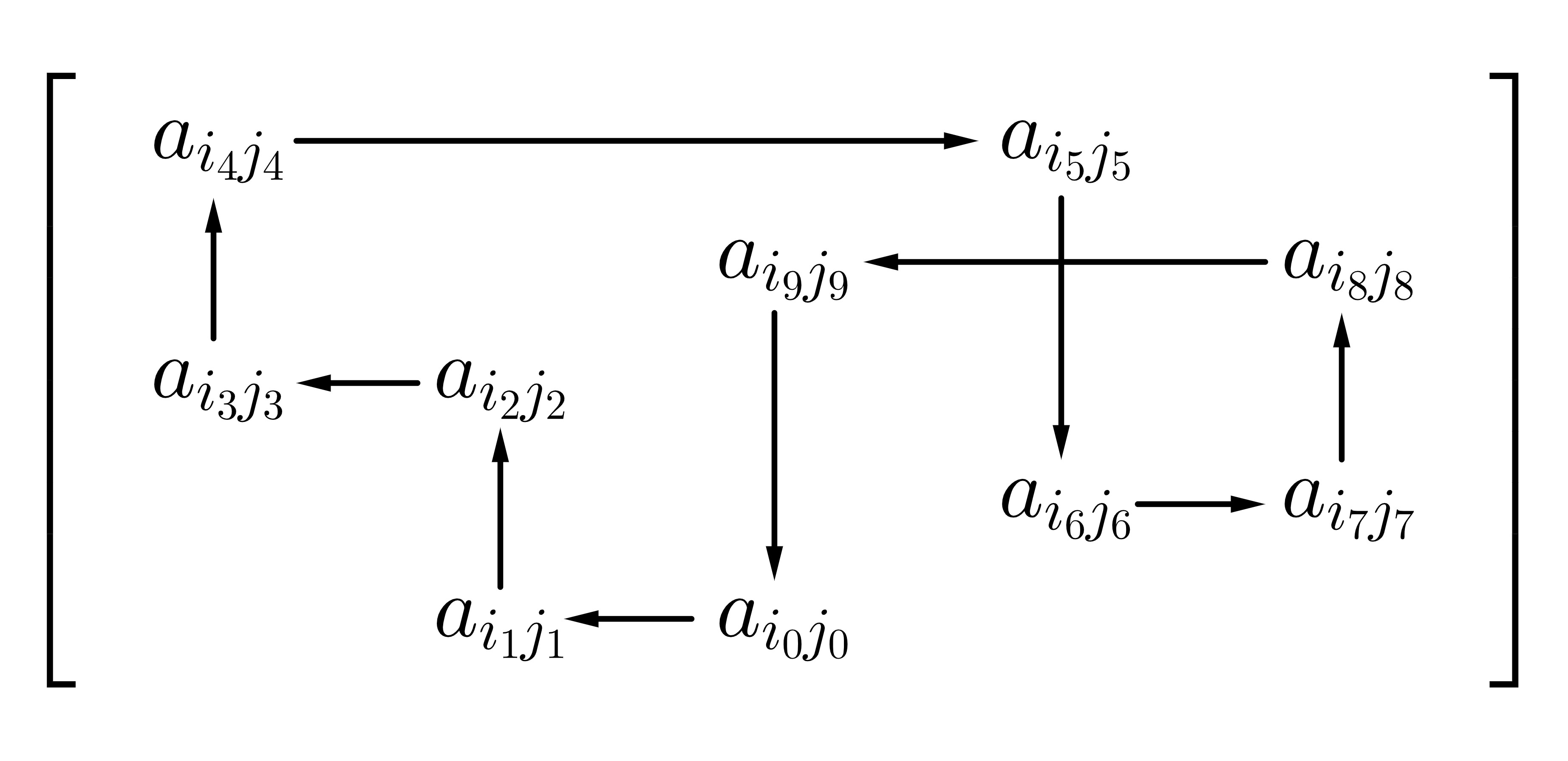}\quad\quad \includegraphics[width=5.4cm]{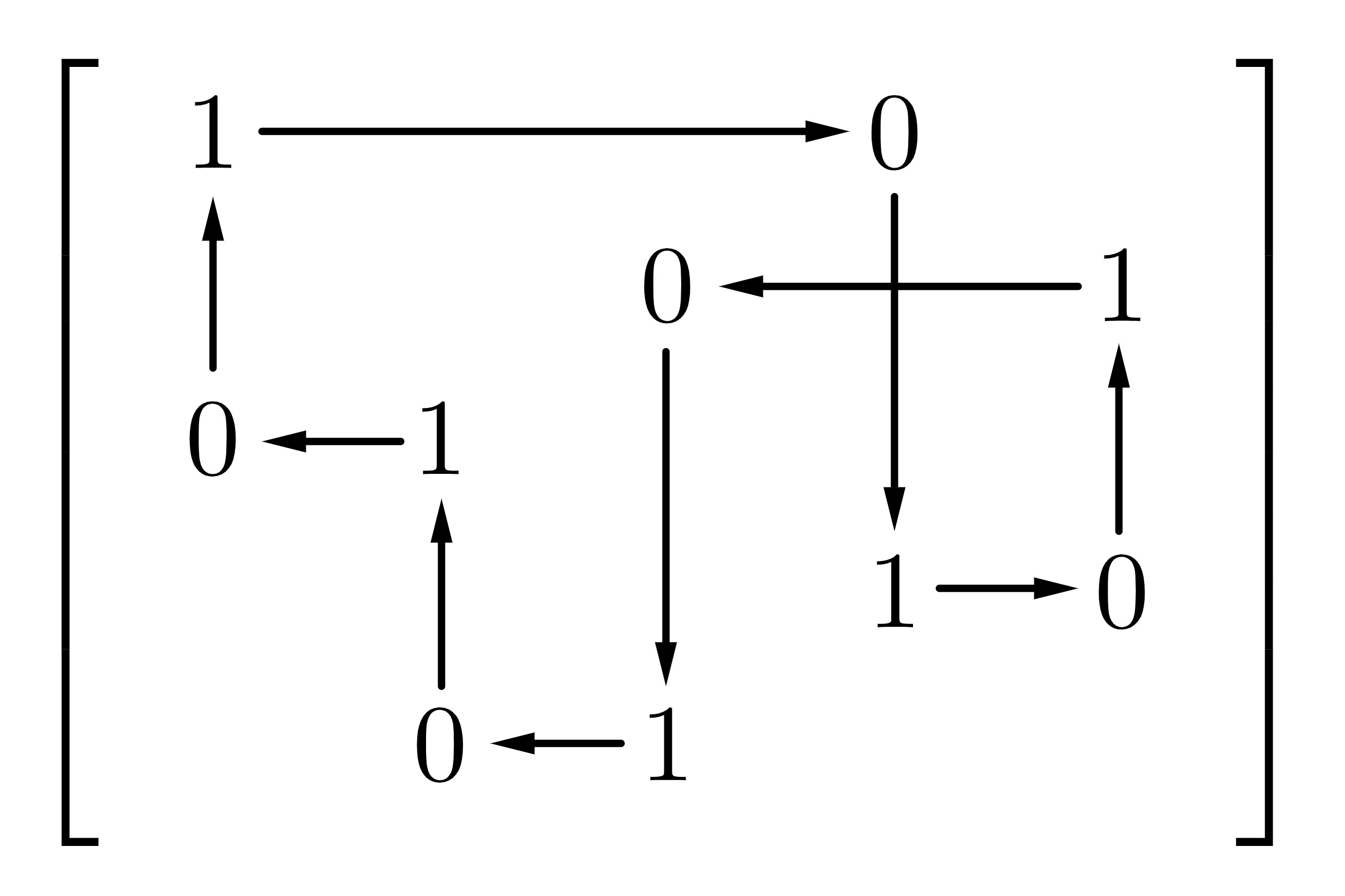}\vspace*{-2mm}
\caption{A switching chain of 10 elements.}
\end{figure}

\noindent see Figure 6. We can find a shortest flow-augmenting path, and hence a shortest switching chain with the above properties by assigning labels to the entries of $A$ based on a breadth-first search in the associated network. This means that we first assign the label $0$ to the entry $a_{ij}$. If we assume that the highest label assigned to any element of $A$ is the even number $k$, then we look for unlabeled entries equal to $0$ in rows containing entries with label $k$. If we find such an entry, then we assign the label $k+1$ to it. If we assume that the highest label assigned to any element of $A$ is the odd number $k$, then we look for unlabeled entries equal to $1$ in columns containing entries with label $k$. If we find such an entry, then we assign the label $k+1$ to it. We repeat these steps as long as possible or the column of $a_{ij}$ contains an entry with the highest non-zero label. Otherwise there's no switching chain that satisfies the above conditions. Let $a_{i_lj_l}=0$ be an entry with the highest label $l$ in the column of $a_{ij}$. Then there must be at least one entry, $a_{i_{l-1}j_{l-1}}=1$, with label $l-1$ in the row of $a_{i_lj_l}$. The column of $a_{i_{l-1}j_{l-1}}$ must contain at least one entry, $a_{i_{l-2}j_{l-2}}=0$, with label $l-2$. This can be continued until we find an entry $a_{i_1,j_1}=0$ with label $1$ in the column of $a_{i_2,j_2}$ and in the row of $a_{ij}$. Thus $(i,j)=(i_0,j_0),(i_1,j_1),\ldots (i_l,j_l)$ is a shortest switching chain that satisfies the above conditions.

It's easy to see, that by interchanging zeros and ones in a switching chain doesn't change the row sums and column sums of the matrix $A$. Hence, if the $i$-th row sum of $A$ is less than the prescribed row sum $r_i$ and the $j$-th column sum of $A$ is less than the prescribed column sum $s_j$, but there exists a switching chain containing the entry $a_{ij}=1$, then interchanging the zeros and ones in the switching chain makes $a_{ij}=0$, and we can switch this to $a_{ij}=1$ to increase the $i$-th row sum and $j$-th column sum of $A$ by $1$. Note that the switch leaves other row sums and column sums unchanged. It's a well-known result in the theory of network flows, that a flow is maximal if and only if there exists no flow-augmenting path. This ensures that, if the tomographic problem has a solution, then a switching chain exists. Therefore, starting with any initial matrix we can find a solution with the help of finitely many switches. The existence of the switching chain can be proved directly with the help of Mirsky's theorem on integer matrices \cite{Mirsky} as the last section shows.

\section{The existence of the switching chain}

Let $R=(r_1,r_2,\ldots,r_m)$ and $S=(s_1,s_2,\ldots,s_n)$ be a pair of compatible integral vectors and let $\mathfrak{A}(R,S)$ denote the set of all $0-1$ matrices of size $m\times n$ with row sums equal to $R$ and column sums equal to $S$. We would like show that if both $\mathfrak{A}(R,S)$ and $\mathfrak{A}(R+\delta,S+\varepsilon)$ are non-empty for some integral vectors $\delta=(\delta_1,\delta_2,\ldots,\delta_m)$ and $\varepsilon=(\varepsilon_1,\varepsilon_2,\ldots,\varepsilon_n)$ such that $R+\delta$ and $S+\varepsilon$ are compatible too, then having a binary matrix $A\in \mathfrak{A}(R,S)$ it's either possible to switch a zero entry of $A$ to 1 without exceeding the row sums in $R+\delta$ and column sums in $S+\varepsilon$, or there exist at least one switching chain in $A$. The proof is based on the following theorem.

\begin{Thm}{Mirsky \cite{Mirsky}}
Let $0\leq r_i'\leq r_i''$, $0\leq s_j'\leq s_j''$ and $c_{ij}\geq 0$ be integers $(i=1,2,\ldots,m)$, $(j=1,2,\ldots,n)$. Then there exists an integral matrix $A=(a_{ij})$ of size $m\times n$ such that
\begin{align*}
r_i'&\leq\sum_{j=1}^na_{ij}\leq r_i''& (i=1,2,\ldots,m)\\
s_j'&\leq\sum_{i=1}^ma_{ij}\leq s_j''& (j=1,2,\ldots,n)\\
0&\leq a_{ij}\leq c_{ij}& (1\leq i\leq m,1\leq j\leq n)
\end{align*}
if and only if, for all $I\subset\{1,2,\ldots,m\}$, $J\subset\{1,2,\ldots,n\}$,
	\[\sum_{i\in I}\sum_{j\in J}c_{ij}\geq\max\left\{\sum_{i\in I}r_i'-\sum_{j\notin J}s_j'',\,\sum_{j\in J}s_j'-\sum_{i\notin I}r_i'' \right\}.
\]
\end{Thm}

\begin{Cor}
The set $\mathfrak{A}(R,S)$ is non-empty with compatible integral vectors $R=(r_1,r_2,\ldots,r_m)$ and $S=(s_1,s_2,\ldots,s_n)$ if and only if for all $I\subset\{1,2,\ldots,m\}$, $J\subset\{1,2,\ldots,n\}$,
\begin{equation}\label{nonempty}
	\left|I\right|\cdot\left|J\right|\geq\sum_{i\in I}r_i-\sum_{j\notin J}s_j=\sum_{j\in J}s_j-\sum_{i\notin I}r_i.
\end{equation}
\end{Cor}

\begin{proof}
Let's choose $r_i'=r_i''=r_i$, $s_i'=s_i''=s_i$, and $c_{ij}=1$ for all $i\in\left\{1,2,\ldots,m\right\}$ and $j\in\left\{1,2,\ldots,n\right\}$. Then, by Mirsky's theorem, the set $\mathfrak{A}(R,S)$ is non-empty if and only if for all $I\subset\{1,2,\ldots,m\}$, and $J\subset\{1,2,\ldots,n\}$,
	\[\left|I\right|\cdot\left|J\right|\geq\max\left\{\sum_{i\in I}r_i-\sum_{j\notin J}s_j,\,\sum_{j\in J}s_j-\sum_{i\notin I}r_i \right\}.
\]On the other hand
\begin{multline*}
	\sum_{i\in I}r_i-\sum_{j\notin J}s_j=\sum_{i\in I}r_i+\sum_{i\notin I}r_i-\sum_{i\notin I}r_i-\sum_{j\notin J}s_j=\\
	\sum_{i=1}^m r_i-\sum_{i\notin I}r_i-\sum_{j\notin J}s_j=\sum_{j=1}^n s_j-\sum_{i\notin I}r_i-\sum_{j\notin J}s_j=\\
	\sum_{j\in J}s_j+\sum_{j\notin J}s_j-\sum_{i\notin I}r_i-\sum_{j\notin J}s_j=\sum_{j\in J}s_j-\sum_{i\notin I}r_i
\end{multline*}
and we are done.
\end{proof}

Given the compatible integral vectors $R=(r_1,r_2,\ldots,r_m)$ and $S=(s_1,s_2,\ldots,s_n)$ let $\delta=(\delta_1,\delta_2,\ldots,\delta_m)$ and $\varepsilon=(\varepsilon_1,\varepsilon_2,\ldots,\varepsilon_n)$ be a pair of integral vectors with $\delta_i\in[0,n-r_i]$ for all $i\in\{1,2,\ldots,m\}$ and $\varepsilon_j\in[0,m-s_j]$ for all $j\in\{1,2,\ldots,n\}$. Then we define the sets $I_0=\{i\,|\,\delta_i>0\}$ and $J_0=\{j\,|\,s_j>0\}$. Let's assume that $I_0$ and $J_0$ are nonempty.

\begin{Thm}\label{deltaeps}
If the sets $\mathfrak{A}(R,S)$ and $\mathfrak{A}(R+\delta,S+\varepsilon)$ are non-empty, then there are indices $i'\in I_0$ and $j'\in J_0$ and there is a matrix $A=\left(a_{ij}\right)$ in $\mathfrak{A}(R,S)$, such that $a_{i'j'}=0$.
\end{Thm}

\begin{proof}
Let's choose arbitrary elements $i_0\in I_0$ and $j_0\in J_0$. By Mirsky's theorem and Corollary 2, there exists a matrix $A\in \mathfrak{A}(R,S)$ with $a_{i_0j_0}=0$ if and only if for all subsets $I\subset\{1,2,\ldots,m\}$, $J\subset\{1,2,\ldots,n\}$,
	\[\sum_{i\in I}\sum_{j\in J}c_{ij}\geq\sum_{i\in I}r_i-\sum_{j\notin J}s_j=\sum_{j\in J}s_j-\sum_{i\notin I}r_i,
\]where
	\[c_{ij}=\left\{
\begin{array}{rl}
0&\textrm{if } i=i_0 \textrm{ and } j=j_0\\
1&\textrm{otherwise}
\end{array}
\right.
\]This means that, if we assume that there's no matrix $A\in \mathfrak{A}(R,S)$ with $a_{i_0j_0}=0$, then there are subsets  $I\subset\{1,2,\ldots,m\}$, $J\subset\{1,2,\ldots,n\}$, such that
\begin{equation}\label{nozero}
\sum_{i\in I}\sum_{j\in J}c_{ij}<\sum_{i\in I}r_i-\sum_{j\notin J}s_j=\sum_{j\in J}s_j-\sum_{i\notin I}r_i.
\end{equation}
It's not possible that $i_0\notin I$ or $j_0\notin J$, since otherwise inequality (\ref{nozero}) means
	\[|I|\cdot |J|<\sum_{i\in I}r_i-\sum_{j\notin J}s_j=\sum_{j\in J}s_j-\sum_{i\notin I}r_i,
\]which implies, that $\mathfrak{A}(R,S)$ is empty. Thus inequality (\ref{nozero}) has the following form
	\[|I|\cdot |J|-1<\sum_{i\in I}r_i-\sum_{j\notin J}s_j=\sum_{j\in J}s_j-\sum_{i\notin I}r_i.
\]This inequality is equivalent to the equation
\begin{equation}\label{nozero2}
|I|\cdot |J|=\sum_{i\in I}r_i-\sum_{j\notin J}s_j=\sum_{j\in J}s_j-\sum_{i\notin I}r_i,
\end{equation}
because $\mathfrak{A}(R,S)$ is non-empty and inequality (\ref{nonempty}) holds. It's not possible that $I_0\subset I$, since otherwise $\sum_{j\in J}\varepsilon_j>0$ and $\sum_{i\notin I}\delta_i=0$, and hence equation (\ref{nozero2}) leads to
	\[|I|\cdot |J|=\sum_{j\in J}s_j-\sum_{i\notin I}r_i<\sum_{j\in J}s_j+\varepsilon_j-\sum_{i\notin I}r_i=\sum_{j\in J}s_j+\varepsilon_j-\sum_{i\notin I}r_i+\delta_i,
\]which implies, that $\mathfrak{A}(R+\delta,S+\varepsilon)$ is empty. Similarly, it's not possible that $J_0\subset J$, since otherwise $\sum_{i\in I}\delta_i>0$ and $\sum_{j\notin J}\varepsilon_j=0$, and hence equation (\ref{nozero2}) leads to
	\[|I|\cdot |J|=\sum_{i\in I}r_i-\sum_{j\notin J}s_j<\sum_{i\in I}r_i+\delta_i-\sum_{j\notin J}s_j=\sum_{i\in I}r_i+\delta_i-\sum_{j\notin J}s_j+\varepsilon_j,
\]which implies, that $\mathfrak{A}(R+\delta,S+\varepsilon)$ is empty. Thus there are elements $i_1\in I_0\setminus I$ and $j_1\in J_0\setminus J$. On the other hand if $A=(a_{ij})$ is any matrix $A\in\mathfrak{A}(R,S)$, then
\begin{multline*}
	\sum_{i\notin I}\sum_{j\notin J}a_{ij}=\sum_{i\notin I}r_i-\sum_{i\notin I}\sum_{j\in J}a_{ij}=\sum_{i\notin I}r_i-\sum_{j\in J}s_j+\sum_{i\in I}\sum_{j\in J}a_{ij}=\\
	\sum_{i\notin I}r_i-\sum_{j\in J}s_j+|I|\cdot |J|-|I|\cdot |J|+\sum_{i\in I}\sum_{j\in J}a_{ij}.
\end{multline*}Hence
	\[\sum_{i\notin I}\sum_{j\notin J}a_{ij}+\left(|I|\cdot |J|-\sum_{i\in I}\sum_{j\in J}a_{ij}\right)=|I|\cdot |J|-\left(\sum_{j\in J}s_j-\sum_{i\notin I}r_i\right).
\]This gives, by equality (\ref{nozero2}), that
	\[\sum_{i\notin I}\sum_{j\notin J}a_{ij}+\left(|I|\cdot |J|-\sum_{i\in I}\sum_{j\in J}a_{ij}\right)=0
\]and here
	\[|I|\cdot |J|-\sum_{i\in I}\sum_{j\in J}a_{ij}\geq 0,
\]thus
	\[\sum_{i\notin I}\sum_{j\notin J}a_{ij}=0.
\]This is possible only if $a_{ij}=0$ for all $i\notin I$ and $j\notin J$, including $i_1\in I_0\setminus I$ and $j_1\in J_0\setminus J$, which gives $a_{i_1j_1}=0$.
Thus finally we can conclude the following.

If $i_0\in I_0$ and $j_0\in J_0$ and there's no matrix $A\in\mathfrak{A}(R,S)$ with $a_{i_0j_0}=0$, then there are indices $i_1\in I_0$ and $j_1\in J_0$, such that $a_{i_1j_1}=0$ for any matrix $A\in \mathfrak{A}(R,S)$. Hence $(i',j')=(i_0,j_0)$ or $(i',j')=(i_1,j_1)$ makes the statement true.
\end{proof}

\medskip
Consider now two compatible integral vectors $R$ and $S$. Assume that none of the row sums of a binary matrix $A=(a_{ij})$ are larger than the corresponding row sums in $R$, and none of the column sums of $A$ are larger than the corresponding column sums in $S$. Let $I_0$ denote the set of all those indexes $i$, where the $i$-th row sum of $A$ is strictly less than the prescribed row sum $r_i$, and let $J_0$ denote the set of all those indexes $j$, where the $j$-th column sum of $A$ is strictly less than the prescribed column sum $s_j$. If $a_{ij}=1$ for all pair of indexes $(i,j)\in I_0\times J_0$, then we can't switch any of the zero entries of $A$ to $1$ without exceeding the prescribed row and column sums. By Theorem \ref{deltaeps}, there must be another binary matrix $\widetilde{A}=(\widetilde{a}_{ij})$ with the same row sums and column sums as $A$, and with $\widetilde{a}_{i'j'}=0$ for some $i'\in I_0$ and $j'\in J_0$ provided that the set $\mathfrak{A}(R,S)$ is nonempty. Ryser showed in \cite{Rys1} that if two matrices have the same row and column sums, then they can be transformed into each other with the help of finitely many switches in so-called switching components, i.e. switching chains of $4$ elements. The entries of $A$ and $\widetilde{A}$ are different in the intersection of $i'$-th row and $j'$-th column, hence $(i',j')$ must be contained in one of the switching components that transform $A$ to $\widetilde{A}$. Merging all these switching components results in a switching chain containing $(i',j')$. Thus Theorem \ref{deltaeps} ensures the existence of the switching chain.

\nocite{*}
\bibliographystyle{fundam}
\bibliography{citations}

%%%%%%%%%%%%%%%%%%%%%%%%%%%%%%%%%%%%%%%%%%%%%%%%%%%%%%%%%%%%%%%%%%%%%%

\end{document}